\documentclass[preprint, 12pt]{article}
\usepackage[colorlinks=true, linkcolor=blue]{hyperref}
\usepackage{graphicx}
\usepackage{amsmath}
\usepackage{amssymb}
\usepackage{mathdots}
\usepackage{amsthm}
\usepackage{soul}
\usepackage{hyperref}
\usepackage[margin=2.5cm]{geometry}
\usepackage{tabularx}
\usepackage{upgreek}
\usepackage{caption}
\usepackage{subcaption}
\usepackage{float}
\usepackage{multirow}
\usepackage{longtable}
\usepackage[table]{xcolor}
\usepackage{booktabs}
\usepackage{enumitem}

\theoremstyle{plain}
\newtheorem{thm}{Theorem}[section]
\newtheorem{lem}[thm]{Lemma}

\newtheoremstyle{myremark}
  {\topsep}   
  {\topsep}   
  {\normalfont}  
  {}         
  {\bfseries} 
  {.}        
  {.5em}     
  {}         

\theoremstyle{myremark}
\newtheorem{rem}{Remark}

\theoremstyle{definition}

\newtheorem{conj}[thm]{Conjecture}

\renewcommand{\epsilon}{\varepsilon}

\makeatletter
\renewcommand{\maketitle}{
  \begingroup
  \renewcommand{\@maketitle}{
    \begin{center}
      {\Large\bfseries\@title\par}
      \vskip 0.8em
      {\normalsize\@author\par}
      \vskip 1.5em
    \end{center}
  }
  \@maketitle
  \thispagestyle{plain}
  \endgroup
}
\makeatother
\begin{document}
\title{Tilings with Infinite Local Complexity and $n$-Fold Rotational Symmetry, $n\in\{13,17,21\}$}

\author{\large{April Lynne D. Say-awen}\\}
\maketitle

\begin{abstract}
A tiling $\mathcal{T}$ is said to have infinite local complexity (ILC) if it contains infinitely many two-tile patches up to rigid motions. In this work, we provide examples of substitution rules that generate tilings with ILC. The proof relies on Danzer’s algorithm \cite{Danz02}, which assumes that the substitution factor is non-Pisot. In addition to ILC, the tiling space of each substitution rule contains a tiling that exhibits (global) $n-$fold rotational symmetry,  $n \in \{13,17,21\}$. 
\end{abstract}
\section{Introduction}

\noindent A \textit{Pisot} (or Pisot--Vijayaraghavan) number is a real algebraic integer \(\nu > 1\) whose conjugates all lie within the unit circle in the complex plane. For example, the golden mean \(\varphi = (1+\sqrt{5})/2\) is a Pisot number because its conjugate \(\varphi' = (1-\sqrt{5})/2 \approx -0.618\) has modulus less than \(1\). Pisot numbers play a significant role in the study of substitution tilings. Specifically, if the substitution factor of a tiling (from a fairly general class of substitution tilings with fractagonal tiles) is a Pisot number, then the tiling exhibits \textit{finite local complexity} (FLC). This result was proven by Frank and Robinson~\cite{FrankRob05}, with a special case given by Frettl\"oh~\cite{Fre02}. A tiling is said to have FLC if it contains only finitely many distinct two-tile patches, up to rigid motions. (Note that some authors define FLC considering only translations.) A tiling that does not exhibit FLC, meaning it has infinitely many two-tile patches, is said to have \textit{infinite local complexity} (ILC).

The converse of the statement “Pisot substitution factor implies FLC” is not generally true. There exist substitution tilings with non-Pisot substitution factors, yet still exhibit FLC. Examples include the chiral Lan\c{c}on--Billard tiling~\cite{LanBill88, BaaGri13} and certain octagonal tilings recently found by Say-awen and Coates~\cite{SayCo}. However, when the additional condition is imposed that some adjacent tiles do not meet full-edge to full-edge, the tiling can be shown to have ILC~\cite{Danz02, FrankRob05}.

In this work, we contribute to the study of substitution tilings by presenting three families of substitution tilings with ILC. We utilize Danzer’s algorithm~\cite{Danz02}, which applies when the substitution factor is non-Pisot and the tiling contains a vertex in a \textit{misfit situation}. A vertex \(\mathcal{V}\) of a tiling \(\mathcal{T}\) is in a misfit situation if there exists a tile \(T \in \mathcal{T}\) such that \(\mathcal{V}\) lies in \(T\) but is not one of its vertices. Equivalently, \(\mathcal{T}\) contains two adjacent tiles that do not meet full-edge to full-edge.

The tilings in this work are constructed using substitution rules whose prototiles are isosceles triangles with interior angles that are multiples of \(\pi/n\), and whose substitution factor is the length of the longest diagonal of a regular \(n\)-gon with unit side length. A more general construction of tilings with this substitution factor and a broader class of prototiles was given by Nischke and Danzer~\cite{NisDanz96}. Similar constructions were presented by Pautze~\cite{Pau17} for \(n = 7, 11\) and by Hibma~\cite{Hib15} for odd \(n\) with \(5 \leq n \leq 15\).

In addition to exhibiting ILC, a distinctive feature of each substitution described here is that its tiling space contains a tiling with \(n\)-fold rotational symmetry. This symmetry distinguishes my examples from most of those in~\cite{NisDanz96} and~\cite{Hib15}. To the best of my knowledge, these are the first examples of substitution tilings with \(n\)-fold rotational symmetry for these specific values of \(n\).

\section{Preliminaries} 
\noindent A \textit{tile} $T \subset \mathbb{R}^2$ is a non-empty compact set that is equal to the closure of its interior. In this work, we restrict our attention to triangular tiles. A \textit{tiling} $\mathcal{T}$ of $\mathbb{R}^2$ is a collection of tiles with pairwise disjoint interiors whose union covers the entire plane. Tiles may be equipped with decorations or markings to indicate their orientations or to distinguish between congruent tiles. Two tiles are \textit{equivalent} if they are congruent and have identical decorations or markings. A finite set $\mathcal{P} \subset \mathcal{T}$ of connected tiles is called a \textit{patch} of $\mathcal{T}$. Two patches of $\mathcal{T}$ are \textit{equivalent} if one can be mapped to the other by an isometry of $\mathbb{R}^2$. For convenience, throughout the paper, whenever we say that a tile $T$ or a patch $\mathcal{P}$ \textit{appears} or \textit{occurs} in $\mathcal{P}'$ (where $\mathcal{P}'$ may be a larger patch, a tiling, or a polygon), we mean that an equivalent copy of $T$ or $\mathcal{P}$ is contained in $\mathcal{P}'$.

An \textit{edge type} of $\mathcal{T}$ is a patch consisting of two tiles that intersect along an edge of $\mathcal{T}$. Recall from the previous section that we defined a tiling $\mathcal{T}$ to exhibit \textit{infinite local complexity} (ILC) if it contains infinitely many two-tile patches. Hence, a tiling containing infinitely many non-equivalent edge types exhibits ILC.

	A basic method for generating a tiling is through a substitution rule. Such a rule consists of a finite set $\mathcal{F}=\{T_1,T_2,…,T_m\}$ of tiles, a real number $\lambda>1$, and a dissection of the inflated tile $\lambda T_{i}$ into tiles $T_{i_1},T_{i_2},…,T_{i_r}$, where each $T_{i_j}$ is equivalent to some tile in $\mathcal{F}$. The mapping $\omega$ that sends each $T_{i}$ to $\{T_{i_1},T_{i_2},…,T_{i_r} \}$ is called \textit{substitution rule} in $\mathbb{R}^2$, with \textit{prototile set} $\mathcal{F}$,  \textit{prototiles} $T_1,T_2,…,T_m$, and \textit{substitution factor} or \textit{inflation factor} $\lambda$. 
	
	The substitution rule $\omega$ naturally extends to any equivalent copy of a prototile. If $T$ is an equivalent copy of the prototile $T_i$ of $\omega$, then $T$ can be written as $T=RT_i+t$, where $R$ is an isometry in the plane that fixes the origin and $t \in \mathbb{R}^2$. We define the image of $T$ under $\omega$ as  $\omega(T)=R\omega(T_i )+ \lambda t$. In a similar way, given $\mathcal{P} \in \mathcal{S}$ --- where $\mathcal{S}$ is the collection of all non-empty sets containing tiles equivalent to the tiles in $\mathcal{F}$ --- $\omega$ extends to any set in $\mathcal{S}$ by $\omega(\mathcal{P})=\{\omega(T)\mid T\in \mathcal{P}\}$. In particular, $\omega$ can be applied $k$ times on $T_i$ to obtain the $k-$\textit{order supertile} $\omega^k (T_i)$ of $T_i$. 

\begin{figure}[H]
\center
\includegraphics[width=.75\textwidth]{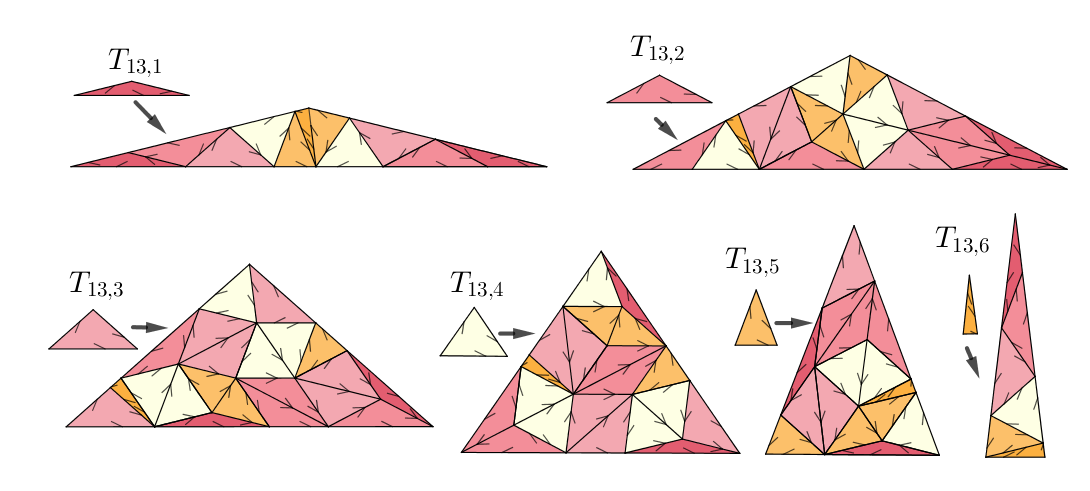}
\caption{The substitution $\sigma_{13}$.  
\label{fig:13-fold}}
\end{figure}

\begin{figure}[H]
\center
\includegraphics[width=1\textwidth]{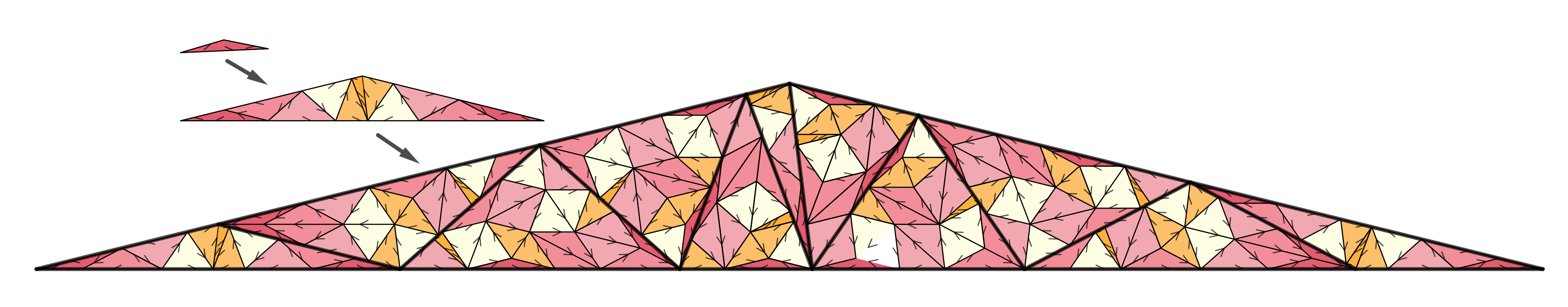}
\caption{The $1$- and $2$-order supertiles of $T_{13,1}$}.  
\label{fig:supertiles}
\end{figure}

As a preliminary illustration, we consider one of the substitution rules which will be examined in detail in a later section. The substitution rule is described in Fig.~\ref{fig:13-fold} with six isosceles triangle prototiles. In particular, the prototile $T_{13,k}$, where $k \in \{1, 2, \ldots, 6\}$, is an isosceles triangle with base angles of $k\pi/13$ and a vertex angle (opening angle) of $(13 - 2k)\pi/13$. The substitution factor is equal to $\mu_{13}=(\sin(6\pi/13))/(\sin(\pi/13)) \approx 4.148$. As shown in the figure, the substitution rules does not preserve the reflection symmetry of the isosceles prototiles, and so it is necessary to equip the tiles with markings to indicate their orientations. Without these markings, one cannot iterate the substitution rule on the prototiles. The 2-order supertile of $T_{13,1}$ is given in Fig.~\ref{fig:supertiles}. 

A patch $\mathcal{P}$ is a \textit{legal patch} of $\omega$ if it appears in some supertile of one of the prototiles of $\omega$. A tiling is a \textit{substitution tiling} with respect to $\omega$ if each of its patches is a legal patch of $\omega$. The collection of all such tilings is called the \textit{tiling space} of $\omega$ and is denoted by $\mathbb{X}_{\omega}$. A common way to obtain a substitution tiling via $\omega$ is to iterate the substitution on a legal patch. In what follows, we illustrate this using the substitution~$\sigma_{13}$.

\begin{figure}[H]
		\begin{subfigure}{1\textwidth}
			\centering
\includegraphics[width=.5\textwidth]{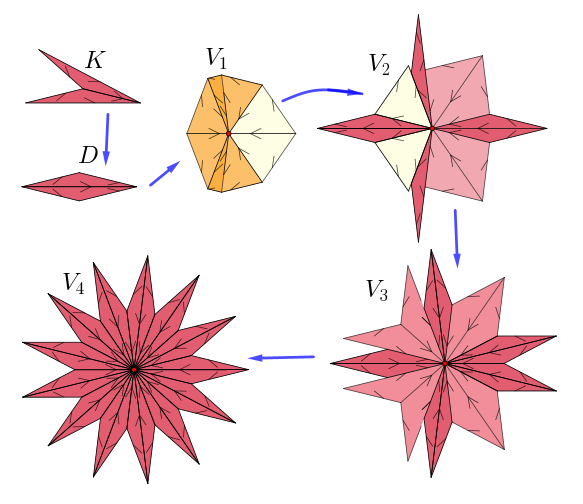}
			\caption{}
		\label{fig:vertexstar}
		\end{subfigure}
		\hfill
		\begin{subfigure}{1\textwidth}
			\centering
			\includegraphics[scale=.6]{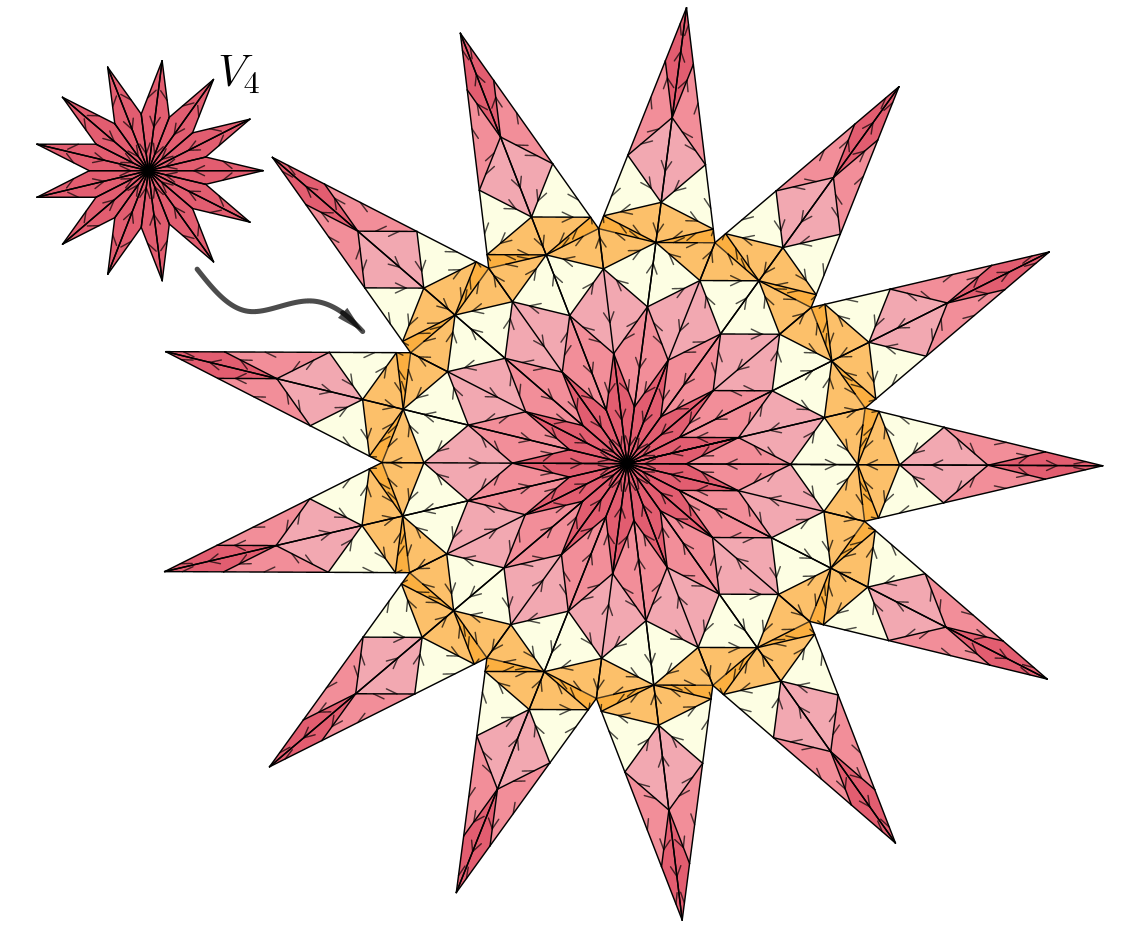}
			\caption{}
		\label{fig:tiling}
		\end{subfigure} 
		\caption{The edge type $K\subset \sigma_{13} (T_{13,2})$ gives rise to the patch $V_4$, and $V_4$ gives rise to the tiling $T_{13}.$}
		\label{fig:rotation}
	\end{figure}

Consider the kite-shaped edge type $K$ located at the right corner of $\sigma_{13}(T_{13,2})$. Applying the substitution to $K$ produces a rhomb-shaped edge type $D$; that is, $\sigma_{13}(K)$ contains a copy of $D$ (see Fig.~\ref{fig:vertexstar}). Applying the substitution to $D$ yields the patch $V_1$. Repeating this process produces the patch $V_4$, which is invariant under $13$-fold rotation. We then center $V_4$ at the origin and apply $\sigma_{13}$ to it, obtaining a larger patch $\sigma_{13}(V_4)$, shown in Fig.~\ref{fig:tiling}. The center of this patch contains a rotated copy of $V_4$; equivalently, $V_4 \subset R \sigma_{13}(V_4),$ where $R$ is the counterclockwise rotation by $\pi/26$ about the origin. Consequently, we obtain the nested sequence
\[
V_4 \subset R \sigma_{13}(V_4) \subset (R\sigma_{13})^2(V_4) \subset (R\sigma_{13})^3(V_4) \subset \dots,
\]

\noindent which converges to a tiling $\mathcal{T}_{13}$. Clearly, $\mathcal{T}_{13} \in \mathbb{X}_{\sigma_{13}}$ since every patch in the nested sequence is legal. Moreover, each patch is invariant under $13$-fold rotational symmetry, and hence so is $\mathcal{T}_{13}$. Finally, $\mathcal{T}_{13}$ is a \textit{fixed point} of $R\sigma_{13}$; that is, $R\sigma_{13}(\mathcal{T}_{13}) = \mathcal{T}_{13}$. \\

Another property of $\sigma_{13}$ is that it is primitive. A substitution rule is said to be \textit{primitive} if, for each prototile $T_i$, there exists a supertile of $T_i$ that contains a copy of every prototile. In the case of $\sigma_{13}$, primitivity is immediate, since each $1$-order supertile already contains at least one copy of each prototile.

An important consequence of the primitivity of a substitution rule $\omega$ is that it guarantees $\mathbb{X}_{\omega}$ is nonempty. Moreover, every legal patch of $\omega$ appears in every tiling $\mathcal{T} \in \mathbb{X}_{\omega}$ (this property is commonly referred to as \textit{local indistinguishability} in the literature; see \cite{BaaGri13}). 

\newpage
\section{The Substitution Rules}
\noindent In the previous section, we introduced one of the substitution rules considered in this work. In what follows, we first present the remaining two substitution rules and the resulting tilings with rotational symmetry. We then describe the Kannan–Soroker–Kenyon (KSK) criterion~\cite{KanSor92, Ken93}, which is used to obtain the dissections that give rise to these substitution rules. We note that substitution rules were not derived through a systematic, step-by-step procedure; rather, they arose from brute-force exploration guided by the KSK criterion and by techniques for enforcing the \(n\)-fold rotational symmetry property. The latter, which constitutes a different line of investigation, follows a distinct approach that is too lengthy to include here and is currently being studied for a general class of tilings in~\cite{SayDe}.  

Our focus here is therefore on describing the resulting constructions and their properties, rather than providing a procedural method for their derivation.

\subsection{The Substitution Rules and Tilings with \(n-\)Fold Rotational Symmetry}
\noindent Recall that the substitution factor for \(n=13\) is  \(\mu_{13} = \frac{\sin\left(\frac{6\pi}{13}\right)}{\sin\left(\frac{\pi}{13}\right)} \approx 4.148,\) which is equal to the length of the longest diagonal of a regular \(13\)-gon with unit side length. In general, the length of the longest diagonal of a regular \(n\)-gon with unit side length is given by  

\begin{equation}
\mu_n = \frac{\sin\left( \frac{(n-1)\pi}{2n} \right)}{\sin\left( \frac{\pi}{n} \right)}, 
\label{eq:subfac}
\end{equation}  
\noindent and we use \(\mu_{17}\) and \(\mu_{21}\) as substitution factors for the cases \(n = 17\) and \(n = 21\), respectively. As mentioned in the introduction, this formula first appeared in \cite{NisDanz96} and has since been used to produce several tilings.

The prototiles are $T_{n,1}, T_{n,2}, \dots, T_{n,(n-1)/2}$, where $T_{n,k}$ is an isosceles triangle with legs of length $1$ and opening angle $(n-2k)\pi/n$. By straightforward calculation, the length of the base of $T_{n,k}$ is
\[
s_{n,k} = 2 \sin\left( \frac{(n-2k)\pi}{2n} \right).
\]

For each $n \in \{13,17,21\}$, we denote the corresponding substitution rule by $\sigma_n$. The case $n = 13$ is already shown in Fig.~\ref{fig:13-fold}, while the substitution rules $\sigma_{17}$ and $\sigma_{21}$ are illustrated in the following figures. As in the case $n = 13$, the tiles are equipped with markings, since the dissections corresponding to the substitution rules break the symmetry of the inflated prototiles, and so markings are necessary to indicate the orientations of the tiles. Moreover, each prototile appears in every 1-order supertile, and so $\sigma_{n}$ is primitive.  

\begin{figure}[H]
\center
\includegraphics[width=.75\textwidth]{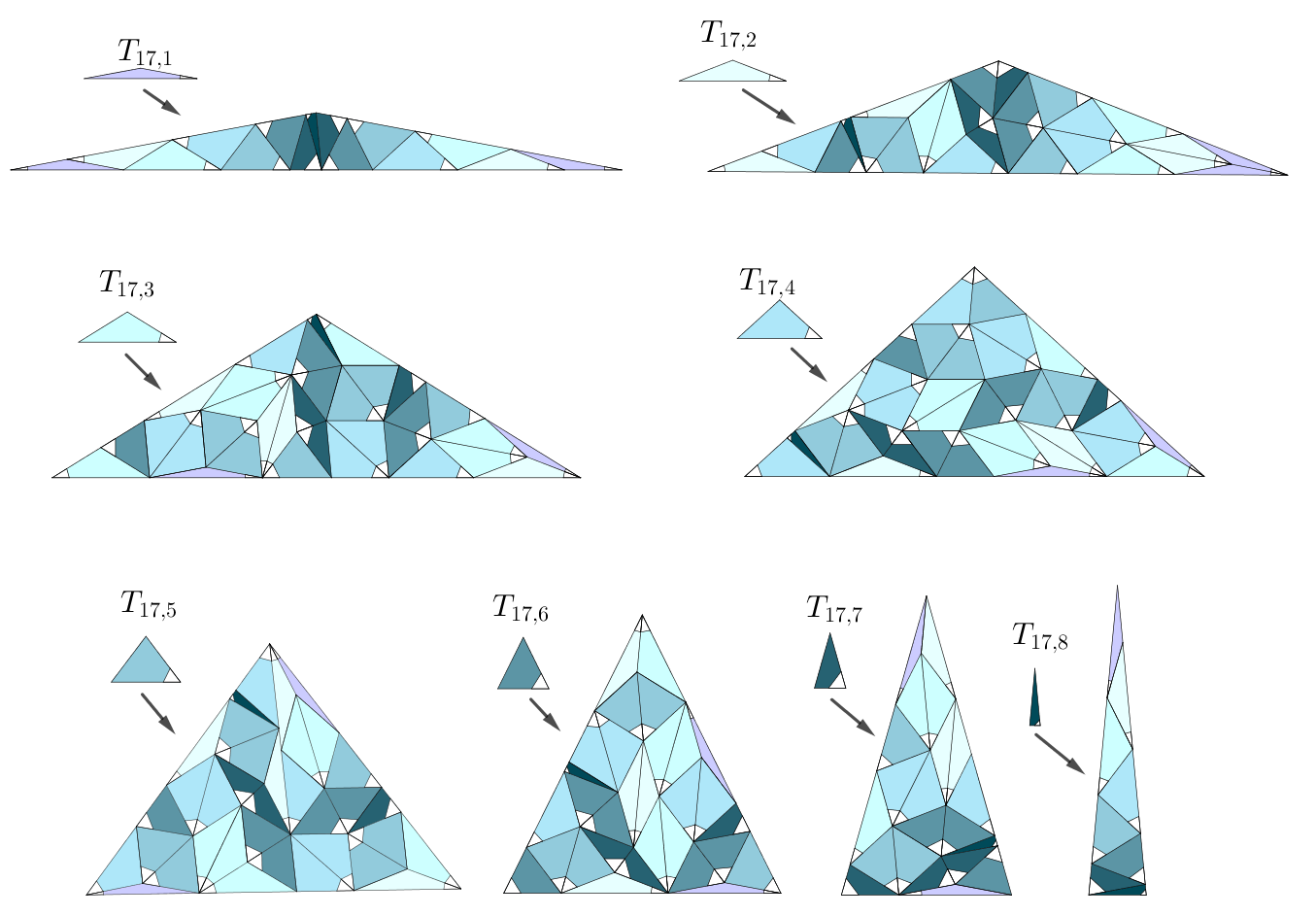}
\caption{The substitution $\sigma_{17}$.  
\label{fig:17-fold}}
\end{figure}

\begin{figure}[H]
\center
\includegraphics[width=.75\textwidth]{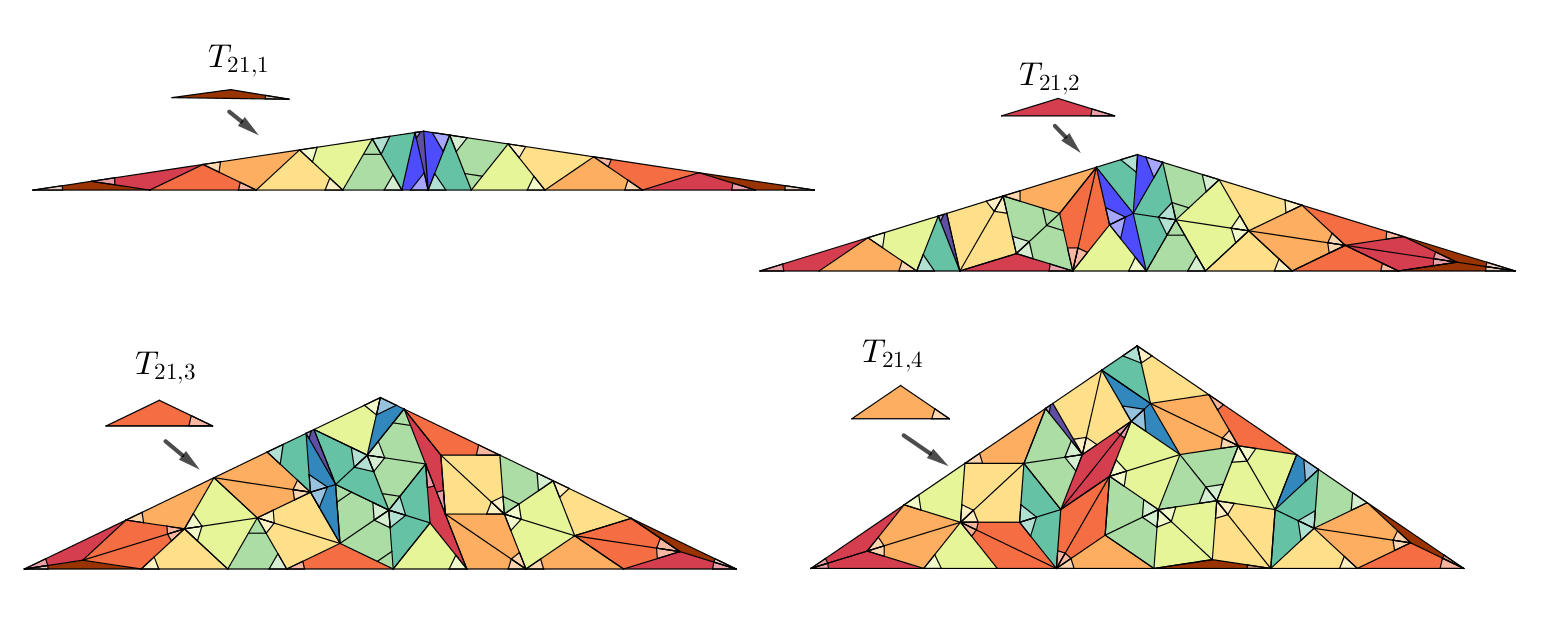}\\
\includegraphics[width=.75\textwidth]{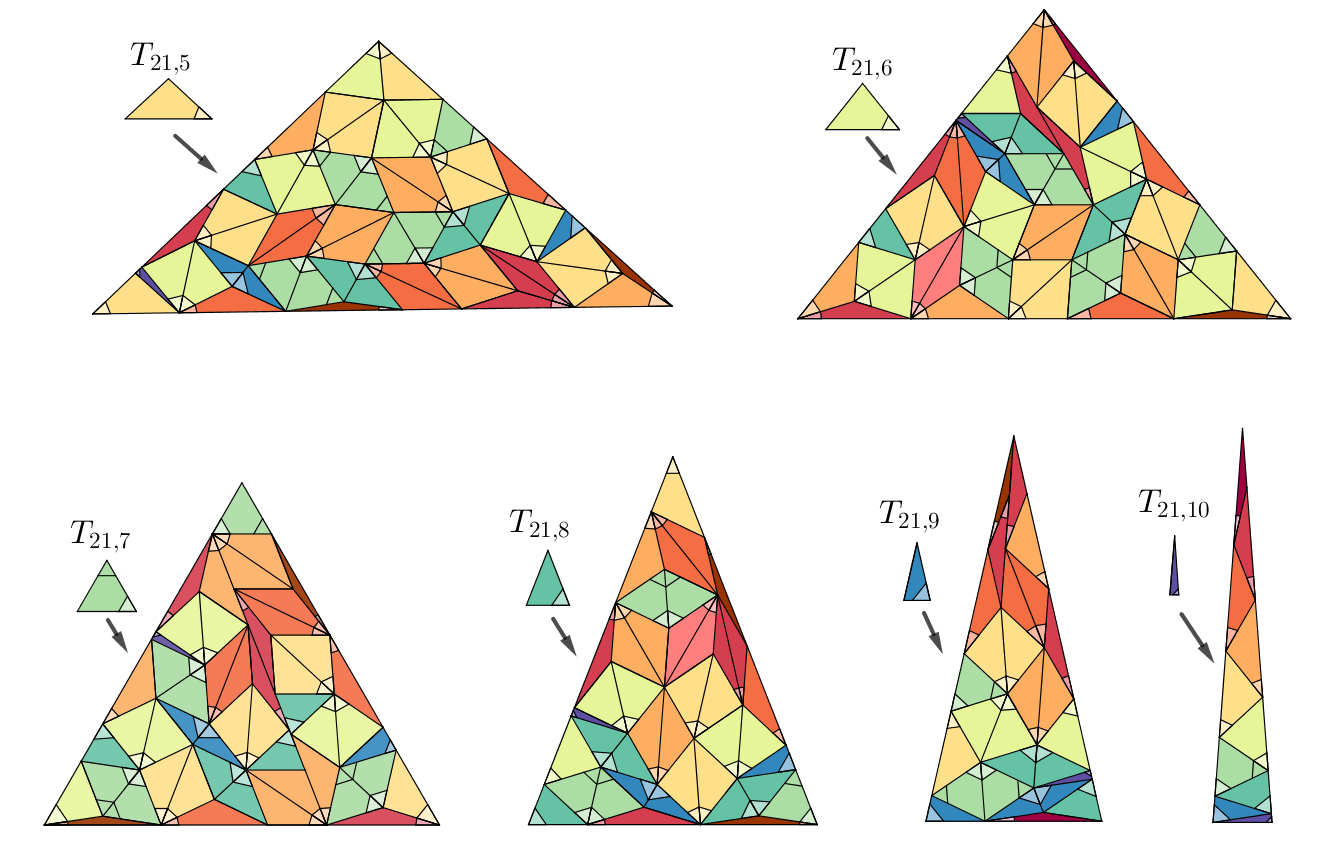}\\
\includegraphics[width=.75\textwidth]{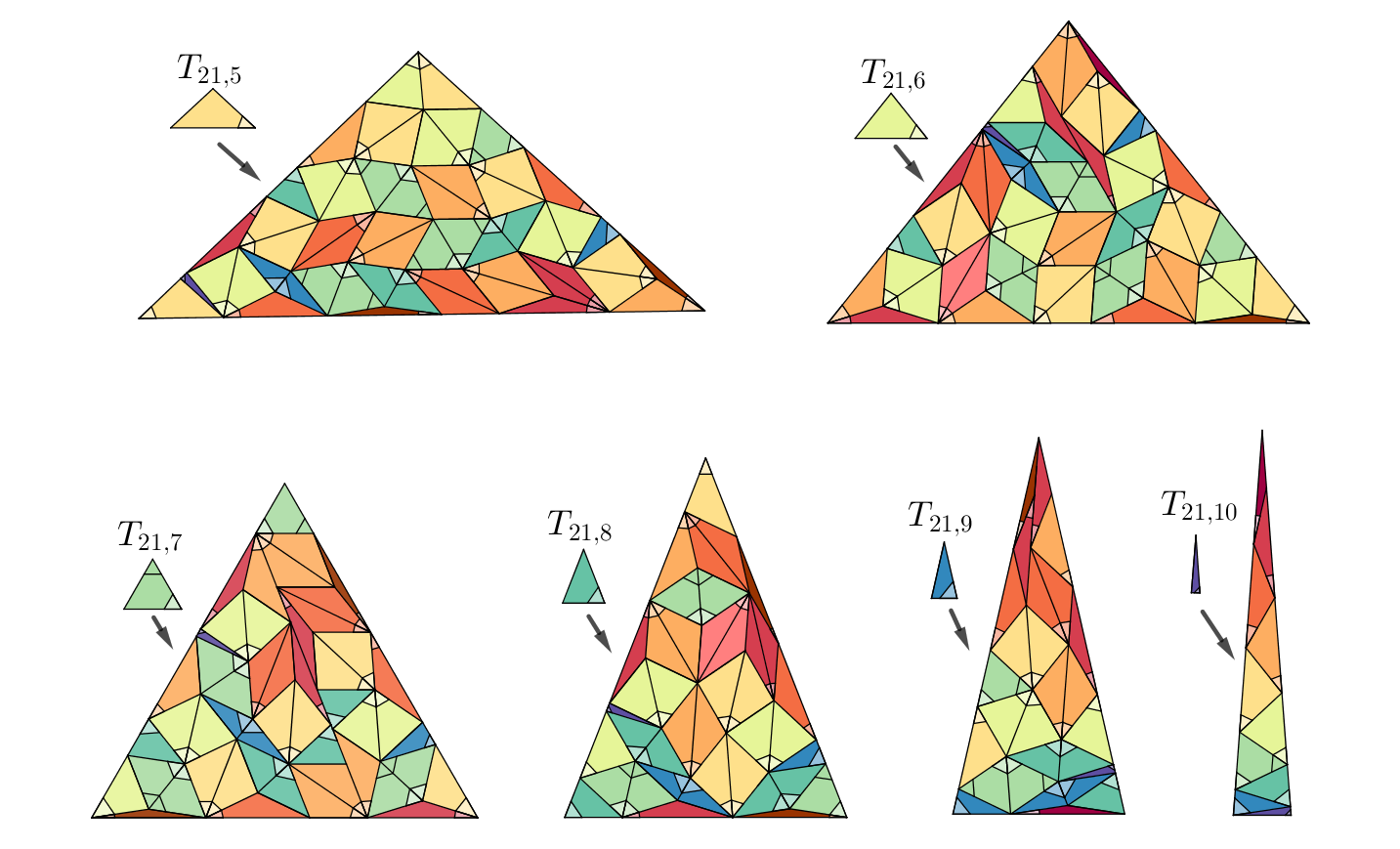}
\caption{The substitution $\sigma_{21}$.  
\label{fig:21-fold}}
\end{figure}
\newpage

As in the case \( n = 13 \), the tiling spaces \(\mathbb{X}_{\sigma_{17}}\) and \(\mathbb{X}_{\sigma_{21}}\) contain tilings \(\mathcal{T}_{17}\) and \(\mathcal{T}_{21}\) with \(17\)-fold and \(21\)-fold rotational symmetry, respectively. The derivation of each tiling is analogous to that of \(\mathcal{T}_{13}\) described in the preliminaries. In each case, the tiling arises from the kite-shaped edge type at the right corner of \(\sigma_{n}(T_{n,2})\), identical to \(K\) in \(\sigma_{13}(T_{13,2})\).

\begin{figure}[H]
\centering
\includegraphics[width=1\textwidth]{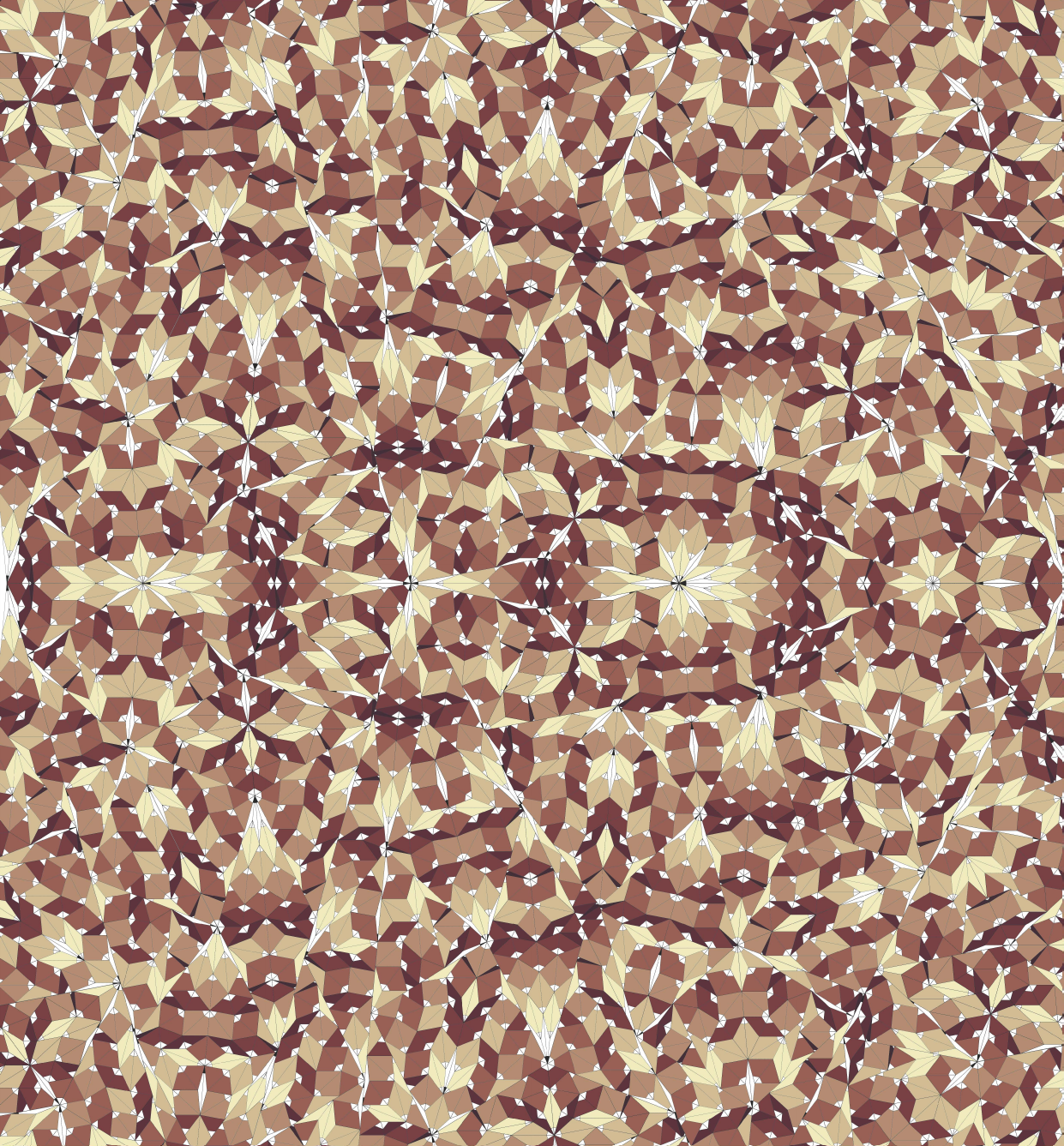}
			\caption{A patch of the tiling $\mathcal{T}_{17}$. The tiles have been recolored for aesthetic purposes.}
\label{fig:17-fold-tiling}
\end{figure}

\begin{figure}[H]
\centering
\includegraphics[width=0.9\textwidth]{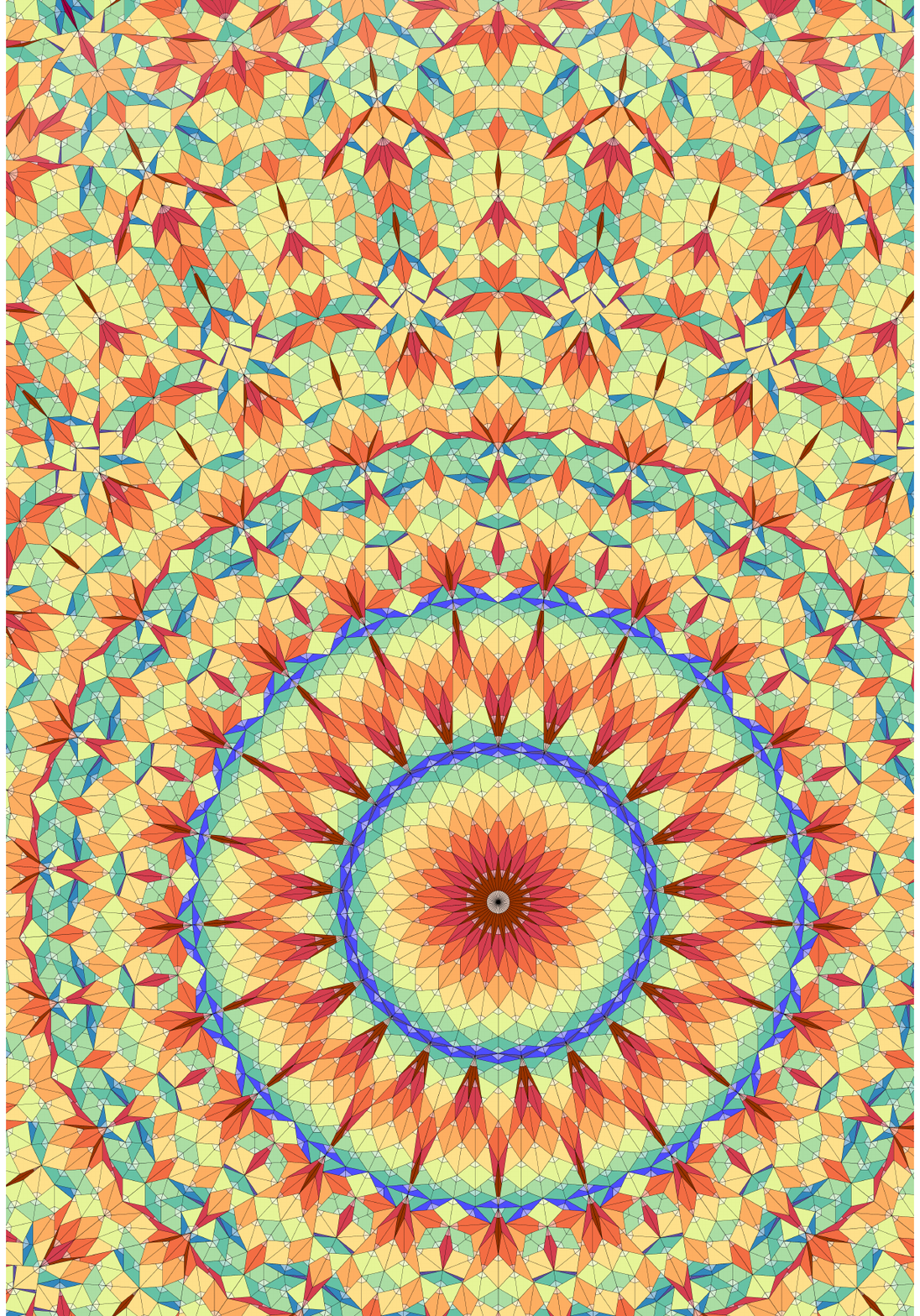}
			\caption{A patch of the tiling $\mathcal{T}_{21}$. }
\label{fig:21-fold-tiling}
\end{figure}

\newpage

\subsection{The KSK Criterion and Its Application to the Construction of $\sigma_n$}
\noindent We now discuss a useful tool---the KSK criterion---which we employed in deriving the substitution rules, in particular in the dissection of the inflated prototiles. This criterion has previously been used to generate substitution tilings with rhombic prototiles (see, for example, \cite{Mal15, KarRis16, Pau17}). Although the prototiles in the present work are triangles, the KSK criterion remains very useful, with the aim of dissecting the inflated prototiles so that isosceles triangles in certain regions of their interior meet along their bases to form rhombs.

The KSK criterion, developed independently by Kannan and Soroker~\cite{KanSor92} and by Kenyon~\cite{Ken93}, determines whether a given polygonal region~$\mathcal{R}$ can be tiled using parallelograms. Here, we assume that the boundary of $\mathcal{R}$ consists of edges of equal length, and hence the resulting parallelograms are rhombs.

\begin{figure}[H]
		\begin{subfigure}{0.45\textwidth}
			\centering
			\includegraphics[scale=.12]{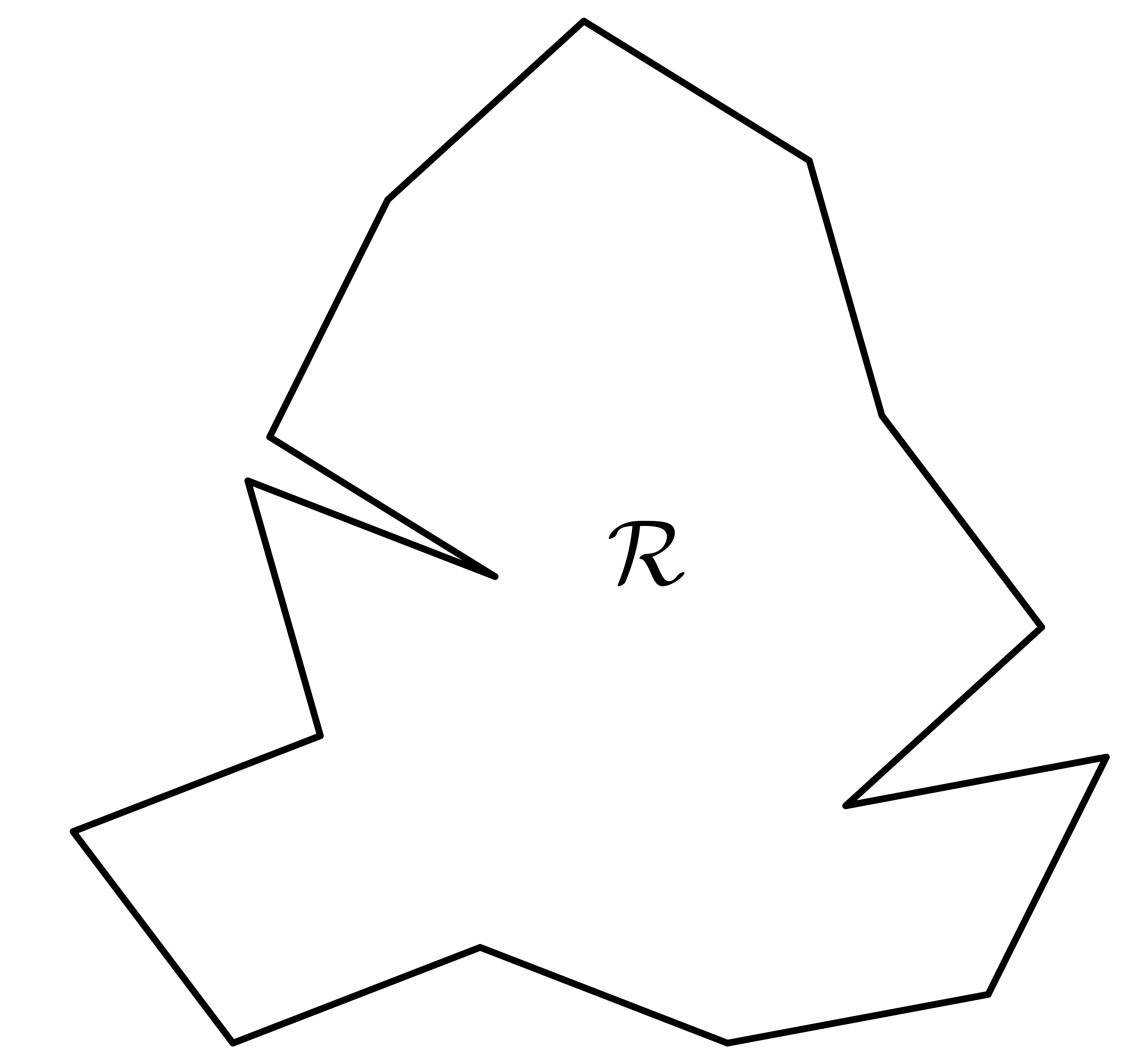}
			\caption{}
\label{fig:convex(a)}
		\end{subfigure}
		\begin{subfigure}{0.45\textwidth}
			\centering
			\includegraphics[scale=.12]{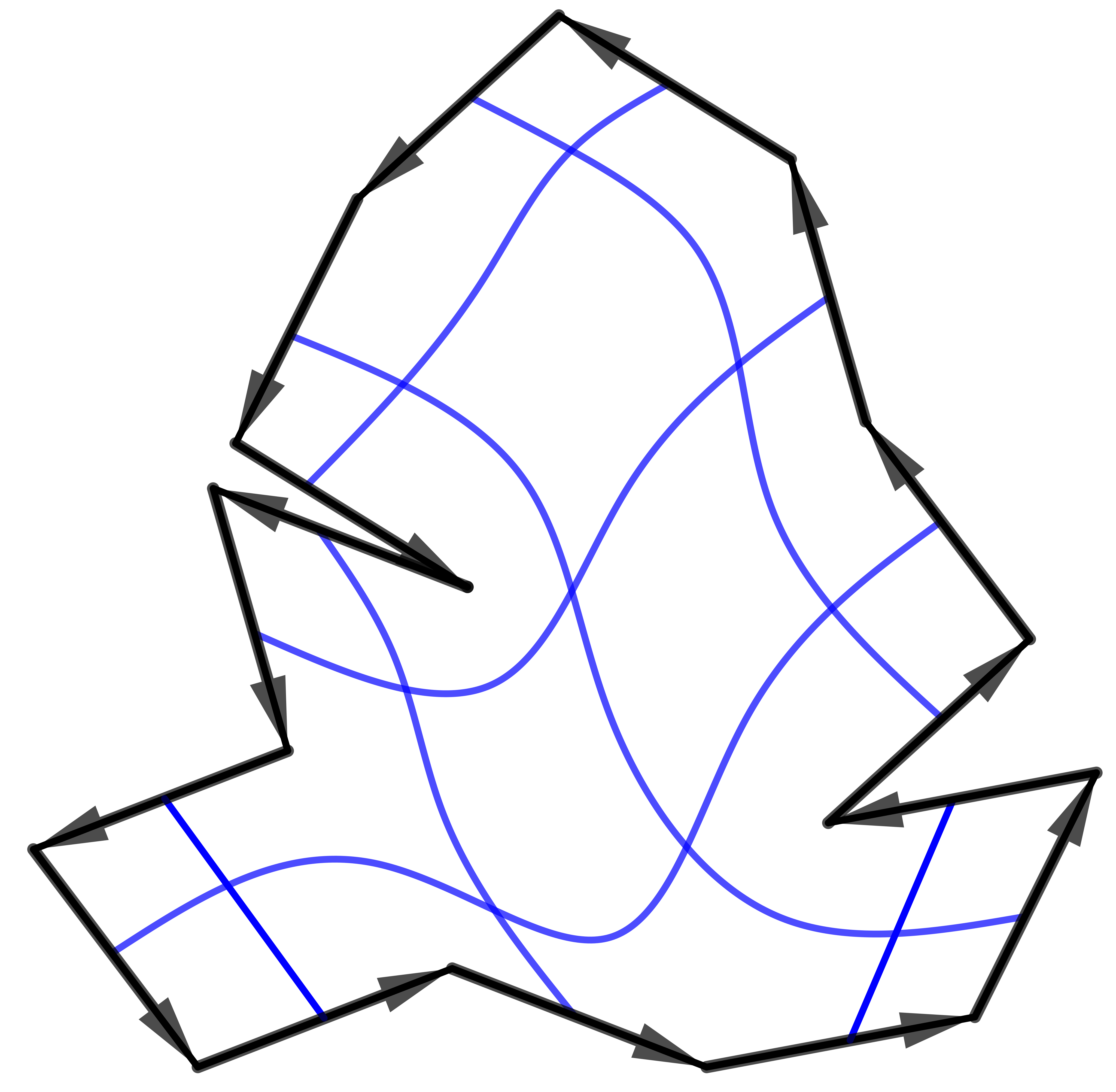}
			\caption{}
\label{fig:convex(b)}
		\end{subfigure} \\
		
		\begin{subfigure}{0.45\textwidth}
			\centering
			\includegraphics[scale=.12]{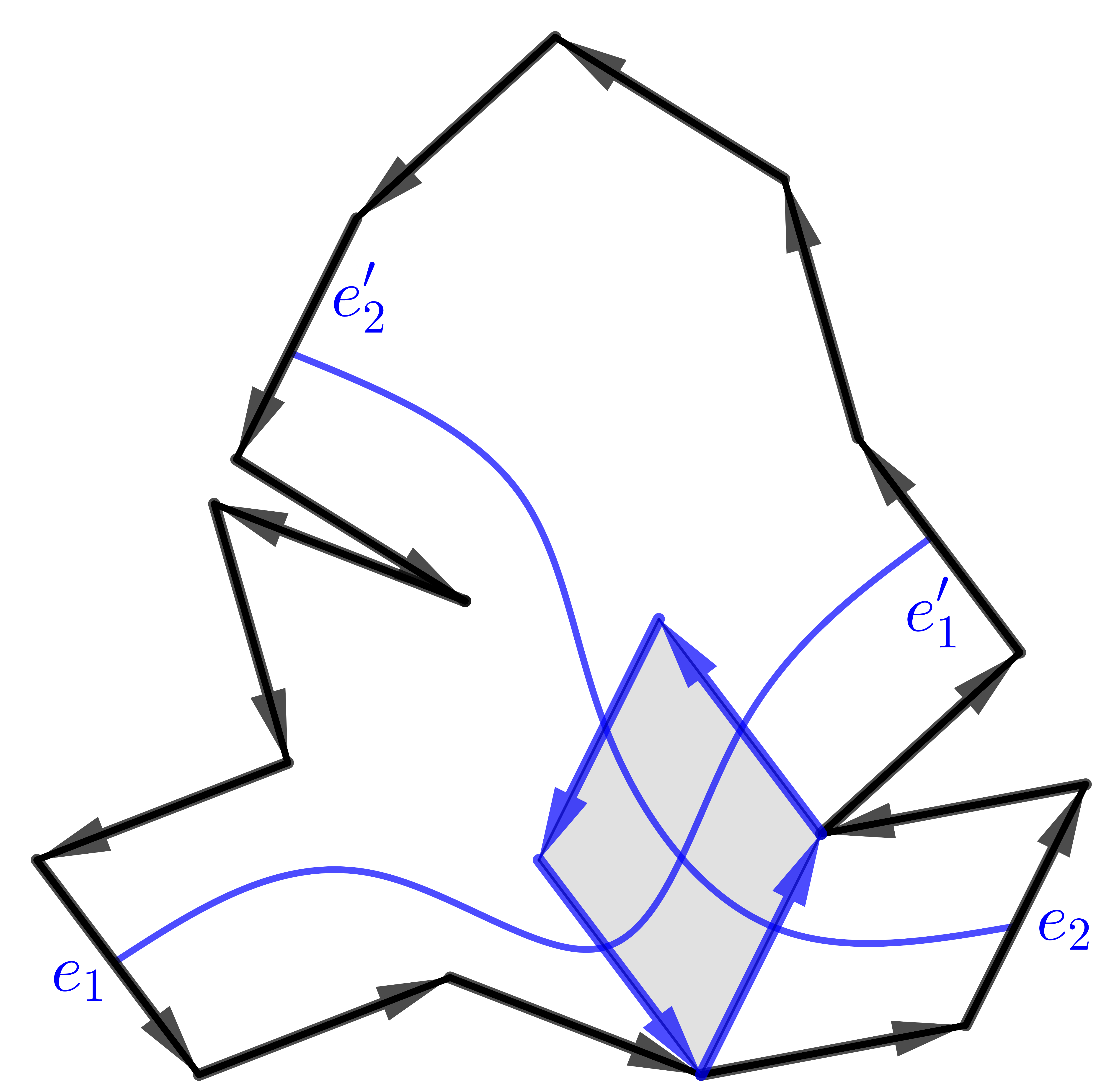}
			\caption{}
\label{fig:convex(c)}
		\end{subfigure}
		\begin{subfigure}{0.45\textwidth}
			\centering
			\includegraphics[scale=.12]{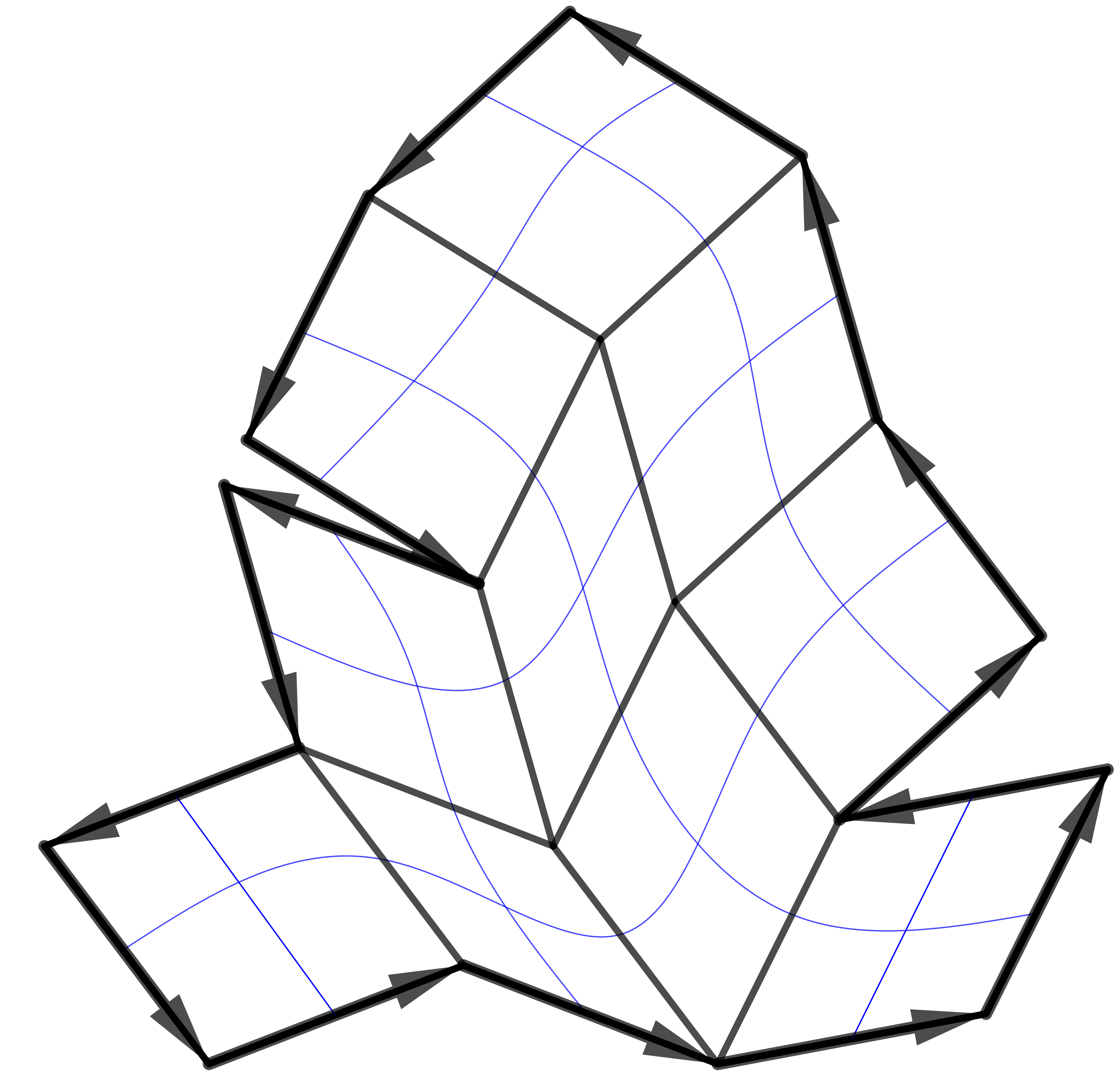}
			\caption{}
\label{fig:convex(d)}
		\end{subfigure}

		\caption{A region that satisfies the KSK criterion.}
		\label{fig:convex}
	\end{figure}

The KSK criterion consists of two conditions, namely, the \textit{balance} and \textit{convex crossing} conditions. We describe the criterion using a concrete example~$\mathcal{R}$ shown in Fig.~\ref{fig:convex(a)}. To properly define these conditions, we view the boundary of~$\mathcal{R}$ as a set~$\mathcal{E}$ of oriented edges arranged in counterclockwise order starting from some fixed vertex (see Fig.~\ref{fig:convex(b)}). The first condition of the criterion is the following. \\

\noindent \textbf{Balance Condition:} The set~$\mathcal{E}$ can be partitioned into pairs of oriented edges with opposite directions.\\

We then join each pair using a smooth planar curve such that each crossing is formed by exactly two curves (see Fig.~\ref{fig:convex(b)}). Our next focus of investigation is on the crossings themselves. Let \(E_1 = \{e_1, e_1'\}\) and \(E_2 = \{e_2, e_2'\}\) be two pairs of edges whose corresponding curves cross, and suppose the sequence of oriented edges in counterclockwise order is given by 
$\ldots e_1 \ldots e_2 \ldots e_1' \ldots e_2' \ldots.$
We say that the pairs \(E_1\) and \(E_2\) have a \textit{convex crossing} if the angles subtended to the left of \(e_1\) and \(e_2\) (when the head of \(e_1\) is attached to the tail of \(e_2\) and the head of \(e_2\) to the tail of \(e_1'\), respectively) are less than \(\pi\). Equivalently, the oriented edges \(e_1, e_2, e_1', e_2'\) form a rhomb inside \(\mathcal{R}\) when the head of \(e_1\) is attached to the tail of \(e_2\), \(e_2\) to \(e_1'\), and \(e_1'\) to \(e_2'\), as illustrated in Fig.~\ref{fig:convex(c)}. With this, the second condition of the criterion is the following:\\

\noindent \textbf{Convex Crossing Condition} Every pair of curves that cross must have a convex crossing.\\

If the two conditions are satisfied, then \(\mathcal{R}\) can be dissected into rhombs. The dissection is carried out by noting that each smooth planar curve corresponds to a chain of rhombs, as illustrated in Fig.~\ref{fig:convex(d)}.

\subsubsection{The Dissection of $\mu_{17}T_{17,6}$}
We note that for any (self‐similar) substitution rule to be well‐defined, the lengths of the inflated sides of the prototiles (obtained by multiplying by the inflation factor) must be expressible as non‐negative integer linear combinations of the prototile edge lengths.

To describe the representations of $\mu_n$ and $\mu_n s_{n,k}$ as sums of $s_{n,j}$, we define
\[
I_n = \{0, 1, 2, \dots, (n-1)/2\}, \quad
I_{n,O} = I_n \cap (2\mathbb{Z}+1), \quad
I_{n,E} = I_n \cap 2\mathbb{Z}.
\]

For each $n \in \{13,17,21\}$, it can be verified that
\begin{equation}
\mu_n = \sum_{j \in I_{n,O}} s_{n,j} \;=\; \sum_{j \in I_{n,E}} s_{n,j},
\label{eq:legrep}
\end{equation}
and
\begin{equation}
\mu_n s_{n,k} = \sum_{j=0}^{\frac{\,n-1\,}{2} - k} s_{n,j}.
\label{eq:baserep}
\end{equation}
These equalities hold for every odd $n \ge 5$ and follow from the identity $(*)$ in~\cite{NisDanz96}.

To illustrate how to obtain a dissection of \(\mu_{n} T_{n,k}\) with the aid of the KSK criterion, we consider the case \(n = 17\) and \(k = 6\).

We define \(T_{17,0}\) as the horizontal segment of length \(s_{17,0} = 1\), and refer to it as a \textit{pseudo-prototile}. An equivalent copy of \(T_{17,0}\) is any isometric image of this segment. Whenever we refer to the \textit{base} or \textit{leg} of \(T_{17,0}\), we simply mean the segment itself.

Using the representations of \(\mu_{17}\) and \(\mu_{17} s_{17,6}\) given in Equations~\ref{eq:legrep} and~\ref{eq:baserep}, each edge of \(\mu_{17} T_{17,6}\) can be dissected into smaller segments, each of length \(s_{17,j}\), where \(j \in I_{17} = \{0, 1, 2, \ldots, 8\}\).  We define the dissection of the left leg using  
$\mu_{17} = s_{17,0} + s_{17,2} + s_{17,4} + s_{17,8},$ and the right leg using $\mu_{17} = s_{17,1} + s_{17,3} + s_{17,5} + s_{17,7}.$ The base is dissected using $s_{17,0} + s_{17,1} + s_{17,2}.$

Noting that \(s_{17,j}\) is the length of the base of the prototile \(T_{17,j}\), we place, for each small edge of length \(s_{17,j}\), a tile equivalent to \(T_{17,j}\) inside the inflated prototile so that its base coincides with the small edge. In this way, any edge dissection \(\mathcal{D}_{17,6}\) determines the tile types along the boundary of \(\mu_{17} T_{17,6}\). Specifically:  
\begin{enumerate}
\item[-] The left leg contains one copy each of \(T_{17,0}, T_{17,2}, T_{17,4}, T_{17,8}\);  
\item[-] The right leg contains one copy each of \(T_{17,1}, T_{17,3}, T_{17,5}, T_{17,7}\);  
\item[-] The base contains one copy each of \(T_{17,0}, T_{17,1}, T_{17,2}\).  
\end{enumerate}

\begin{figure}[H]
  \centering
  \includegraphics[scale=.11]{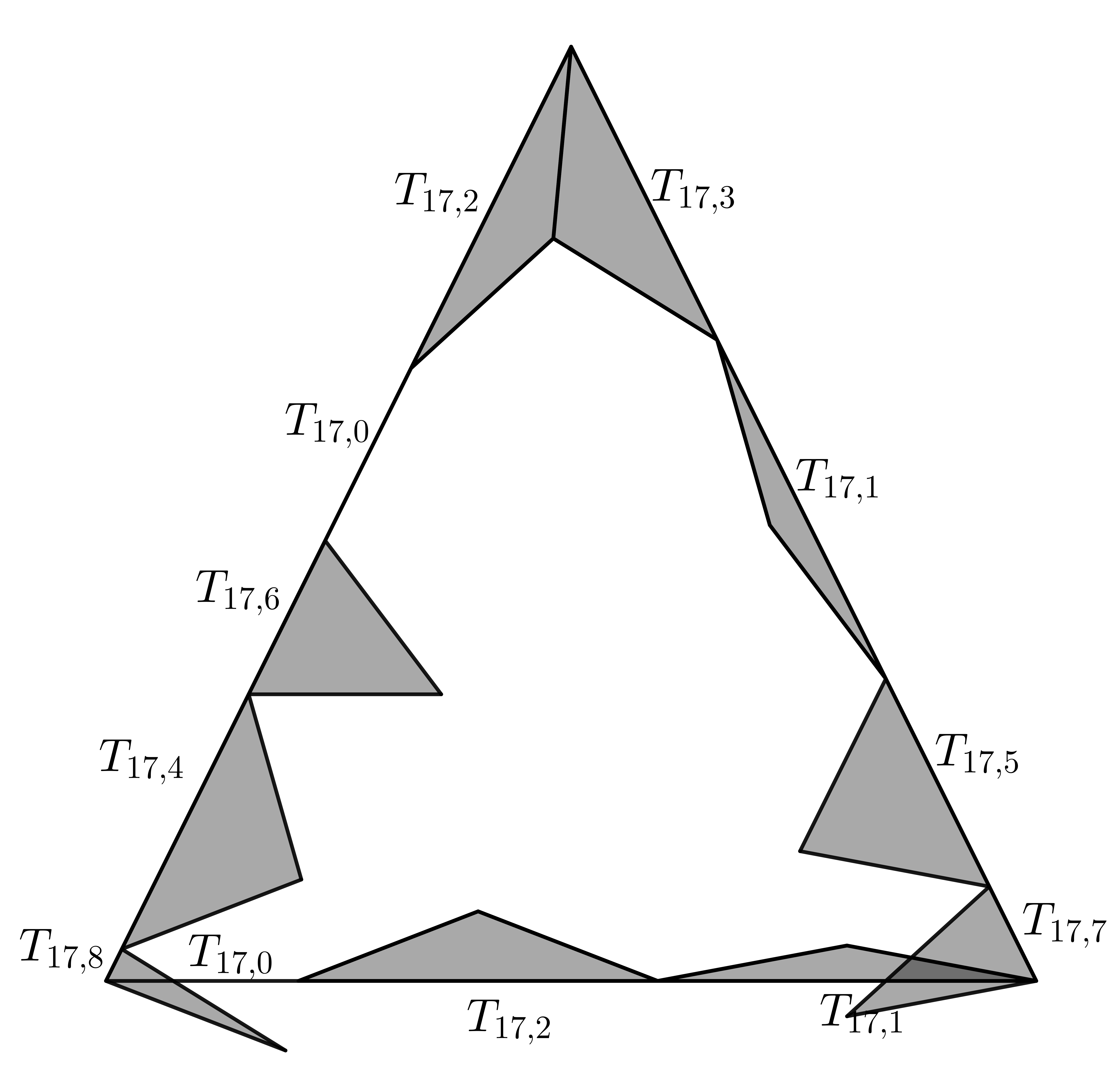}
  \caption{An inadmissible arrangement of tiles along the boundary.}
  \label{fig:wrongarr}
\end{figure}

The question now is: how should \(\mathcal{D}_{17,6}\) be chosen to ensure that the inflated triangle \(\mu_{17} T_{17,6}\) can be dissected into copies of the prototiles? An obvious restriction on \(\mathcal{D}_{17,6}\) is that the initial (boundary) tiles must have pairwise disjoint interiors and must lie entirely within the large triangle \(\mu_{17} T_{17,6}\). For instance, the tile equivalent to \(T_{17,8}\) on the left side of \(\mu_{17} T_{17,6}\) cannot be placed at a corner; otherwise, part of it would lie outside the large triangle (see Fig.~\ref{fig:wrongarr}). Beyond these geometric constraints, the KSK criterion is crucial for identifying admissible choices of \(\mathcal{D}_{17,6}\).

At this stage, \(\mathcal{D}_{17,6}\) is not yet specified. For any given choice of \(\mathcal{D}_{17,6}\), let \(R_{17,6}\) denote the region that remains to be dissected after placing the corresponding initial tiles along the boundary of \(\mu_{17} T_{17,6}\). Our objective is to determine \(\mathcal{D}_{17,6}\) so that the resulting \(R_{17,6}\) satisfies the KSK criterion.

We begin by assigning orientations to the legs of the initial tiles in a counterclockwise direction. We then show that, for any choice of \(\mathcal{D}_{17,6}\), the corresponding region \(R_{17,6}\) satisfies the balance condition. Table~\ref{tab:labels} lists the prototiles that appear along each side of \(\mu_{17} T_{17,6}\), together with the positive angles their legs make with the \(x\)-axis. As shown in the table, all these angles are distinct and can be partitioned into pairs, as detailed in Table~\ref{tab:partition}, where the difference between the angles in each pair is either \(\pi\) or \(-\pi\), indicating opposite directions. For example, the left leg of the copy of \(T_{17,2}\) on the left side of the large triangle is paired with the right leg of the copy of \(T_{17,3}\) on the right side, since $\frac{25\pi}{17} - \frac{8\pi}{17} = \frac{17\pi}{17} = \pi.$

Matched edges may either overlap completely or appear as edges of the boundary of the undissected region \(R_{17,6}\). The oriented edges of the boundary of \(R_{17,6}\) are precisely those that do not overlap, yet can still be grouped into pairs of oriented edges with opposite directions. Hence, the balance condition is satisfied.

\begin{center}
\renewcommand{\arraystretch}{1.5}
\captionof{table}{Prototiles that appear along the sides of $\mu_{17}T_{17,6}$ and the angles $s\pi/17$ made by the left (L) and right (R) legs of the tiles with the positive $x$-axis, measured counterclockwise.}
\label{tab:labels}
\begin{tabular}{| c | l | c | l | c | l |}
\hline
Left & Angles & Right & Angles & Base & Angles \\
& $L, R$ & & $L, R$ & & $L, R$ \\
\hline
$T_{17,0}$ & $\frac{23\pi}{17}$ & $T_{17,1}$ & $\frac{12\pi}{17}, \frac{10\pi}{17}$ & $T_{17,0}$ & $0$ \\
$T_{17,2}$ & $\frac{25\pi}{17}, \frac{21\pi}{17}$ & $T_{17,3}$ & $\frac{14\pi}{17}, \frac{8\pi}{17}$ & $T_{17,1}$ & $\frac{\pi}{17}, \frac{33\pi}{17}$ \\
$T_{17,4}$ & $\frac{27\pi}{17}, \frac{19\pi}{17}$ & $T_{17,5}$ & $\frac{16\pi}{17}, \frac{6\pi}{17}$ & $T_{17,2}$ & $\frac{2\pi}{17}, \frac{32\pi}{17}$ \\
$T_{17,6}$ & $\frac{29\pi}{17}, \pi$ & $T_{17,7}$ & $\frac{18\pi}{17}, \frac{4\pi}{17}$ & & \\
$T_{17,8}$ & $\frac{31\pi}{17}, \frac{15\pi}{17}$ & & & & \\
\hline
\end{tabular}
\end{center}

\begin{center}
\renewcommand{\arraystretch}{1.5}
\captionof{table}{Pairs of angles that differ by $\pm \pi$.}
\label{tab:partition}
\begin{tabular}{| c | c | c | c | c | c | c | c | c | c | c | }
\hline
$0$ & $\frac{2\pi}{17}$ & $\frac{4\pi}{17}$ & $\frac{6\pi}{17}$ & $\frac{8\pi}{17}$ & $\frac{10\pi}{17}$ & $\frac{12\pi}{17}$ & $\frac{14\pi}{17}$ & $\frac{16\pi}{17}$ & $\frac{18\pi}{17}$ & $\frac{32\pi}{17}$ \\ 
$\pi$ & $\frac{19\pi}{17}$ & $\frac{21\pi}{17}$ & $\frac{23\pi}{17}$ & $\frac{25\pi}{17}$ & $\frac{27\pi}{17}$ & $\frac{29\pi}{17}$ & $\frac{31\pi}{17}$ & $\frac{33\pi}{17}$ & $\frac{\pi}{17}$ & $\frac{15\pi}{17}$ \\
\hline
\end{tabular}
\end{center}

\begin{figure}[H]
		\begin{subfigure}{0.45\textwidth}
			\centering
			\includegraphics[scale=1]{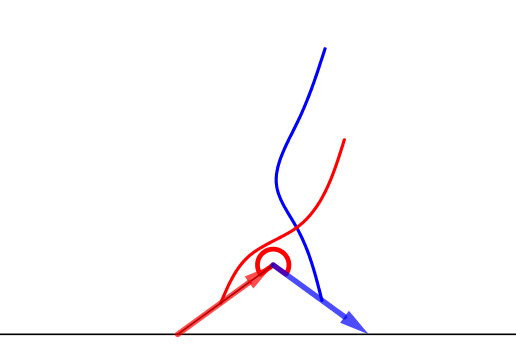}
			\caption{}
\label{fig:diss(a)}
		\end{subfigure}
		\hfill
		\begin{subfigure}{0.45\textwidth}
			\centering
			\includegraphics[scale=1]{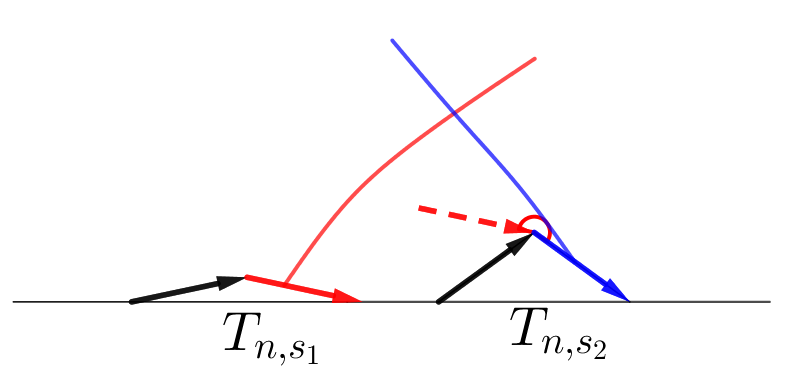}
			\caption{}
\label{fig:diss(b)}
		\end{subfigure} \\
		
		\begin{subfigure}{0.45\textwidth}
			\centering
			\includegraphics[scale=1]{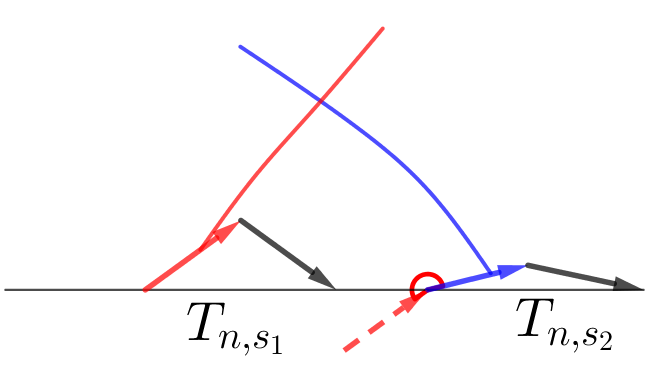}
			\caption{}
\label{fig:dissrec(c)}
		\end{subfigure}
		\hfill
		\begin{subfigure}{0.45\textwidth}
			\centering
			\includegraphics[scale=1]{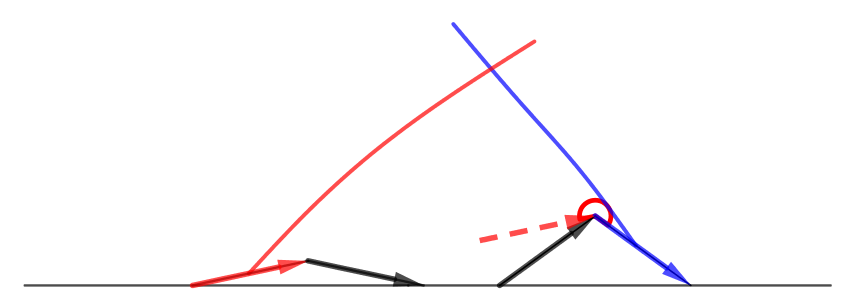}
			\caption{}
\label{fig:dissrec(d)}
		\end{subfigure}
		\caption{Non-convex crossings scenarios. In case (b), $s_{1}<S_{2}$; whereas in case (c), $s_{1}>s_{2}$.}
		\label{fig:nonconvex}
	\end{figure}

Next, we arrange the tiles with the goal of producing a region \(R_{17,6}\) that satisfies the convex-crossing condition. Since the large triangle has only three sides, two of the four oriented edges that define a crossing must lie on the same side. Figure~\ref{fig:nonconvex} illustrates scenarios in which two smooth planar curves, emerging from two edges on the same side of \(\mu_{n} T_{n,k}\), form a non-convex crossing. The corresponding angles---clearly greater than \(\pi\)---are marked in the figure.  

\begin{figure}[H]
			\centering
			\includegraphics[scale=.11]{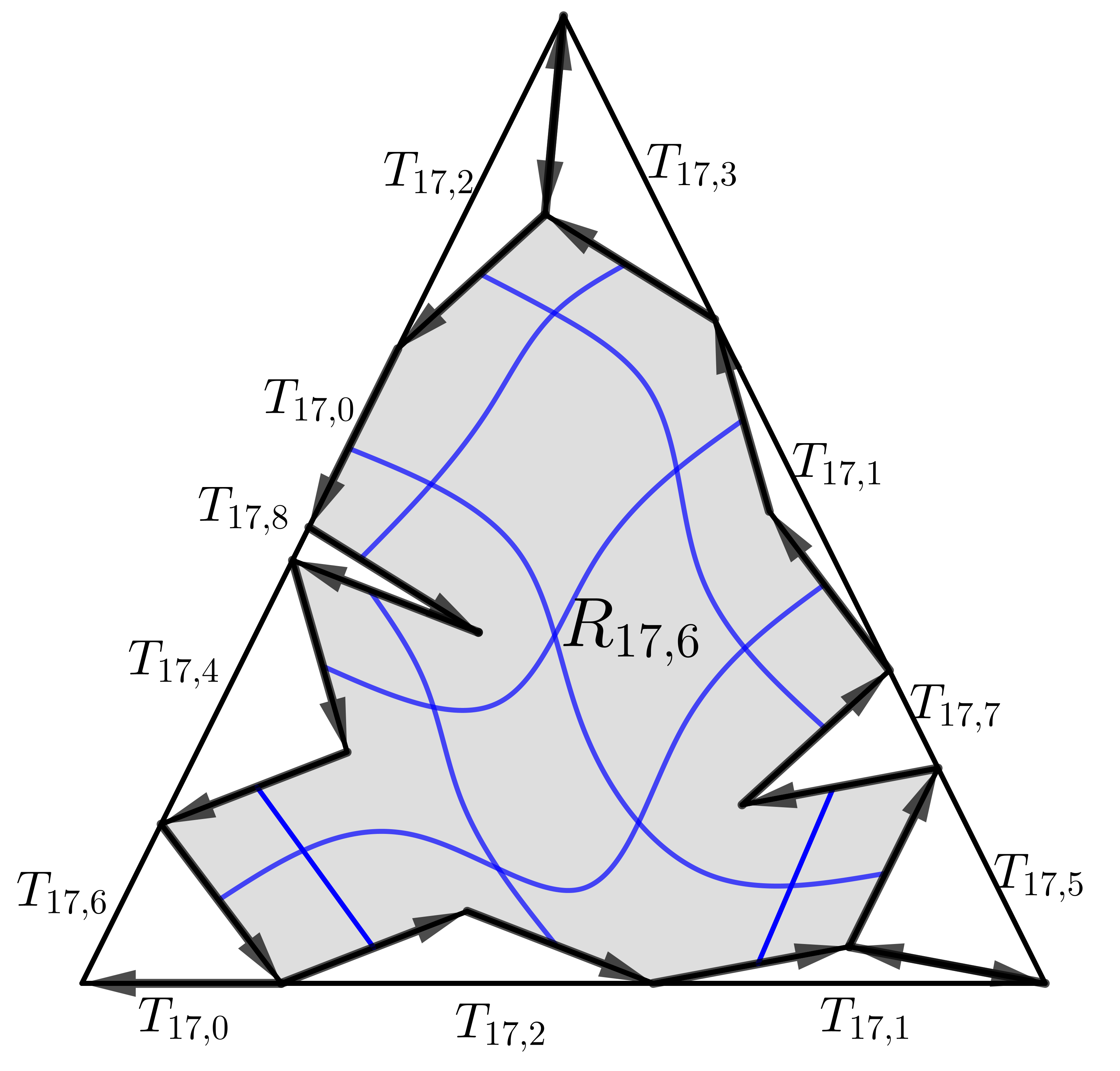}
		\caption{A region $R_{17,6}$ that satisfies the convex crossing.}
		\label{fig:correctdissection}
	\end{figure}

We arrange the tiles from the base of the large triangle to the apex while avoiding the four scenarios. With this, we obtain the arrangement shown in Fig.~\ref{fig:correctdissection}, which contains no convex crossings. The resulting region \(R_{17,6}\) is congruent to \(\mathcal{R}\) in Fig.~\ref{fig:convex}. Consequently, the region can be dissected into rhombs, and each rhomb further subdivided into equilateral triangles congruent to the prototiles.

\begin{rem}
For the cases \(k = 1\) and \(k = (n-1)/2\), the total area of the initial tiles is equal to the area of \(\mu_{n} T_{n,k}\). In this situation, matched edges overlap. Since at least one of the interior angles is \(\pi/n\), the tiles are straightforward to arrange once a tile equivalent to \(T_{n,1}\) is placed at the left or right corner of the large triangle for \(k = 1\), or at the apex for \(k = (n-1)/2\). Using the matching of edges determined by the corresponding angles, the placement of the remaining tiles proceeds in a cascading manner, much like a domino effect.
\end{rem}

\begin{rem}
We note that all values of \(n\) and \(k\) satisfy the balance condition, whereas a few values of \(k\) do not satisfy the convex-crossing condition. We have observed that \(R_{n,k}\) is easier to dissect when the number of non-convex crossings is minimized. In fact, even when a non-convex crossing occurs, it is still possible to associate chains of rhombs with the curves forming the non-convex crossing, although these chains are interrupted along a line (see, for example, \(\sigma_{17}(T_{17,5})\) and \(\sigma_{17}(T_{21,3})\)). This observation suggests further study. Based on examples, including some not presented here, we conjecture that a version of the KSK criterion for isosceles triangles can be formulated.
\end{rem}

\section{Algebraic Properties of $\mu_{n}$} 
\noindent In this section, we discuss algebraic properties of $\mu_n, n \in \{13,17,21\}$, which are essential for showing that a tiling arising from $\sigma_n$ has ILC. 

\begin{lem}
\label{lem:mu_n_properties}
For each $n \in \{13,17,21\}$, the inflation factor $\mu_n$ is a non-Pisot algebraic integer, and the minimal polynomial $g_n(x)$ of $\mu_{n}$ has no repeated roots.
\end{lem}

\begin{proof}
The representations for $\mu_n$ and $\mu_n s_{n,k}$ as sums of prototile edge lengths, given in Equations~\ref{eq:legrep} and~\ref{eq:baserep}, naturally give rise to the matrix equation  $\mathbf{z}^\top M_n = \mu_n \mathbf{z}^\top,$ where 
\[
\mathbf{z} = 
\begin{bmatrix}
s_{n,0}=1 \\ s_{n,1} \\ s_{n,2} \\ \vdots \\
s_{n,(n-3)/2} \\ s_{n,(n-1)/2}
\end{bmatrix} \  \mbox{and} \ 
M_n = 
 \begin{bmatrix} 1 & 1 & 1 & \cdots &1  & 1 & 1 \\ 
0 & 1 & 1 & \cdots & 1 & 1 & 0 \\ 
1 & 1 & 1 & \cdots & 1 & 0 & 0 \\ 
0 & 1 & 1 & \cdots & 0 & 0 & 0 \\ 
\vdots & \vdots & \vdots & \reflectbox{$\ddots$} & \vdots & \vdots & \vdots \\ 
0 & 1 & 0 & \cdots & 0 & 0 & 0 
\\ 1 & 0 & 0 & \cdots & 0 & 0 & 0 \\ \end{bmatrix}.
\]

Each component of the row vector $\mathbf{z}^\top M_n$ corresponds to a representation of $\mu_n s_{n,k}$, $k=0, 1, \ldots, (n-1)/2$, in terms of the edge lengths $s_{n,j}$. For example, the first entry encodes the representation  $s_{n,0} + s_{n,2} + \cdots + s_{n,(n-1)/2}$ of $\mu_n$ as stated in Equation~\ref{eq:legrep}, while the second entry corresponds to  $s_{n,0} + s_{n,1} + \cdots + s_{n,(n-3)/2} = \mu_n s_{n,1}$ following Equation~\ref{eq:baserep}.  

This confirms that $\mu_n$ is an eigenvalue of the integer matrix $M_n$, with left eigenvector $\mathbf{z}^\top$. Hence, $\mu_{n}$ is an algebraic integer, and its minimal polynomial $g_{n}(x)$ divides the characteristic polynomial $p_n(x)$ of $M_n$.

To determine $g_n(x)$ explicitly, we compute the characteristic polynomial  $p_n(x) = \det(M_n - xI)$, where $I$ denotes the identity matrix of size $(n+1)/2$. For each $n =13,17,21$, we have

\begin{equation}
p_n(x)=-\left(x-1 \right) \left(x+\dfrac{1}{s_{n,1}}\right)\left(x-\dfrac{1}{s_{n,2}} \right) \cdots \left(x+\dfrac{1}{s_{n,(n-3)/2}} \right)\left(x-\dfrac{1}{s_{n,(n-1)/2}} \right),
\end{equation}

\noindent with $\mu_{n} = \dfrac{1}{s_{n,(n-1)/2}}$. Thus, the eigenvalues are the positive or negative reciprocals of the prototile edge lengths. Note that the $(n+1)/2$ roots are all distinct, so each eigenvalue is simple. So, $g_{n}(x)$ has no repeated roots. For $n = 13$ and $n = 17$, the minimal polynomial is $g_{n}(x) = \frac{p_{n}(x)}{-(x - 1)},$ while for $n = 21$ we have  $g_{21}(x) = \frac{p_{21}(x)}{-\,(x^2 - 1)\,(x^3 + 2x^2 - x - 1)}.$ These polynomials are listed in Table~\ref{tab:minpoly}, together with the algebraic conjugate $\mu_{n,2}$ of $\mu_n$ whose modulus exceeds $1$, confirming that $\mu_n$ is non-Pisot.
\newpage
\begin{center}
\captionof{table}{The minimal polynomial of $\mu_n$ and its conjugate $\mu_{n,2}$ with $|\mu_{n,2}|>1.$}
\renewcommand{\arraystretch}{1.5}
\label{tab:minpoly}
\begin{tabular}{| c | c | c |}
\hline
$n$ & Minimal Polynomial $g_n (x)$ &	$\mu_{n,2}$ \\
\hline
$13$ & $x^6-3x^5-6x^4+4x^3+5x^2-x-1$ & $-\dfrac{1}{s_{13,5}} \approx -1.41002$ \\ \hline
$17$ & $x^8-4x^7-10x^6+10x^5+15x^4-6x^3-7x^2+x+1$ & $-\dfrac{1}{s_{17,7}} \approx -1.82706$ \\ \hline
$21$ & $x^6-8x^5+8x^4+6x^3-6x^2-x+1$ & $\dfrac{1}{s_{21,8}} \approx 1.36858$ \\
\hline
\end{tabular}
\end{center}
\end{proof}

The non-Pisot property of the inflation factor will become crucial later in the proof of ILC, where it will allow us to construct an infinite sequence of distinct numbers. These numbers will come from lengths of certain segments in the tiling, and our first step toward that goal is to consider the algebraic setting in which these lengths live.

Since $\mu_n$ is an algebraic integer, it generates the ring $\mathbb{Z}[\mu_n]$, which is a free $\mathbb{Z}$-module of rank $d_n = \deg g_n(x)$. 
This module can be naturally identified with $\mathbb{Z}^{d_n}$ via the standard (monomial) basis $B = \{1, \mu_n, \mu_n^2, \dots, \mu_n^{d_n - 1}\}.$

For each $s \in \mathbb{Z}[\mu_n]$, let  
$$
s = a_0 + a_1 \mu_n + a_2 \mu_n^2 + \cdots + a_{d_n-1} \mu_n^{d_n-1}, 
\quad a_i \in \mathbb{Z},
$$  
be its unique representation in terms of the basis $B$. The mapping  
$$
s \longmapsto (a_0, a_1, \dots, a_{d_n-1})^\top
$$  
is the \textit{coordinate isomorphism} with respect to $B$, and $(a_0, a_1, \dots, a_{d_n-1})^\top$ is called the \textit{coordinate vector} of $s$ with respect to $B$. For simplicity, whenever we refer to the coordinate vector of an element of $\mathbb{Z}[\mu_n]$, we mean with respect to the standard basis $B$.

In the next lemma, we establish that $s_{n,j} \in \mathbb{Z}[\mu_n]$ and then compute their coordinate vectors.

\begin{lem}
For each \( n \in \{13, 17, 21\} \) and \( j \in I_n=\{0,1,2,\ldots,(n-1)/2\} \), we have \( s_{n,j} \in \mathbb{Z}[\mu_n] \).
\label{lem:s_n,j}
\end{lem}

\begin{proof}
Let \(\mathbf{o}\) denote the zero vector of dimension \((n+1)/2\). Since \(\mathbf{z}\) is an eigenvector of \(M_n\) with eigenvalue \(\mu_n\) by Lemma~\ref{lem:mu_n_properties}, it satisfies the homogeneous system
\[
(M_n^\top - \mu_n I)\,\mathbf{x} = \mathbf{o}.
\]

Solving this system by Gaussian elimination expresses each \(s_{n,j}\) as a \(\mathbb{Q}\)-linear combination of powers of \(\mu_n\). Reducing these powers using the relation \(g_n(\mu_n) = 0\) shows that each \(s_{n,j}\) can be written in the standard basis \(B\) of \(\mathbb{Z}[\mu_n]\). The corresponding coordinate vectors are listed in Table~\ref{tab:vectorrep}.
\end{proof}

As a final preparatory step,  we recall how multiplication in $\mathbb{Z}[\mu_n]$ by $\mu_n$ corresponds to $\mathbb{Z}^{d_n}$. Let $G_n$ denote the companion matrix of $g_n(x)$. For $s \in \mathbb{Z}[\mu_n]$ with coordinate vector $(a_0, a_1, \dots, a_{d_n - 1})^\top$, the vector $G_n(a_0, a_1, \dots, a_{d_n - 1})^\top$ is precisely the coordinate vector of $\mu_n s$. 

By Lemma~\ref{lem:mu_n_properties}, the polynomial $g_n(x)$ has no repeated roots. Consequently, the companion matrix $G_n$ is diagonalizable over the splitting field of $g_n$. Thus there exists an invertible matrix $A$ such that
\begin{equation}
\label{eq:Gtilde}
A^{-1} G_n A = \widetilde{G}_n =
\begin{bmatrix}
\mu_n & 0 & 0 & \cdots & 0 \\
0 & \mu_{n,2} & 0 & \cdots & 0 \\
0 & 0 & \mu_{n,3} & \cdots & 0 \\
\vdots & \vdots & \vdots & \ddots & \vdots \\
0 & 0 & 0 & \cdots & \mu_{n,d_n}
\end{bmatrix},
\end{equation}
where the $i$th column of $A$ is a right eigenvector and the $i$th row of $A^{-1}$ is a left eigenvector of $G_n$ corresponding to the eigenvalue $\mu_{n,i}$, for $i = 1, 2, \ldots, d_n$, with $\mu_{n,1} = \mu_n$. We also note that $\mu_{n,2}$ is the algebraic conjugate of $\mu_n$ with modulus greater than $1$, as given in Table~\ref{tab:minpoly}.

In this work, we choose $A^{-1}$ so that its second row is
\[
v_{n,2} =
\left( \frac{1}{\mu_{n,2}^{d_n - 1}},\ \frac{1}{\mu_{n,2}^{d_n - 2}},\ \ldots,\ \frac{1}{\mu_{n,2}},\ 1 \right).
\]
As will be needed later, we also include \( v_{n,2} y_{n,j} \) in Table~\ref{tab:vectorrep}.
\newpage

\begin{table}[htbp]
\centering
\caption{The coordinate vector $y_{n,j}$ of $s_{n,j}$ with respect to the standard basis $B$; and $v_{n,2} {y}_{n,j}$. }
\label{tab:vectorrep}
\begin{tabular}{|l|l|l|}
\hline
\( s_{n,j} \) & Coordinate vector $y_{n,j}$ & \( v_{n,2} {y}_{n,j} \) \\
\hline
\( s_{13,0} \) & \( (1,0,0,0,0,0)^\top \) & -0.179421119170724 \\ \hline
\( s_{13,1} \) & \( (-4,1,10,-3,-4,1)^\top \)& -0.268597272296381 \\ \hline
\( s_{13,2} \) & \( (5,-4,-20,11,11,-3)^\top \) & -0.0432536521449611 \\ \hline
\( s_{13,3} \) & \( (0,5,5,-11,-6,2)^\top \) & 0.203845625239822 \\ \hline
\( s_{13,4} \) & \( (-5,0,16,-5,-8,2)^\top \) & 0.348414935064387 \\ \hline
\( s_{13,5} \) & \( (4,-5,-15,14,10,-3)^\top \) & 0.317739022198404 \\ \hline
\( s_{13,6} \) & \( (-1,5,4,-6,-3,1)^\top \) & 0.127247211393141 \\
\hline
\( s_{17,0} \) & \( (1,0,0,0,0,0,0,0)^\top \) & -0.0147137185440312 \\ \hline
\( s_{17,1} \) & \( (-6,1,21,-5,-20,6,5,-1)^\top \) & -0.0250197112728618 \\ \hline
\( s_{17,2} \) & \( (14,-6,-70,29,85,-33,-24,5)^\top \) & -0.0131169374224207 \\ \hline
\( s_{17,3} \) & \( (-14,14,84,-64,-126,67,41,-9)^\top \) & 0.00271522134344693 \\ \hline
\( s_{17,4} \) & \( (0,-14,-14,56,42,-48,-20,5)^\top \) & 0.0177339928490863 \\ \hline
\( s_{17,5} \) & \( (14,0,-70,14,90,-28,-25,5)^\top \) & 0.02744026786708 \\ \hline
\( s_{17,6} \) & \( (-14,14,78,-70,-117,71,40,-9)^\top \) & 0.0289263790503149 \\ \hline
\( s_{17,7} \) & \( (6,-14,-35,55,56,-45,-21,5)^\top \) & 0.0217471384189736 \\ \hline
\( s_{17,8} \) & \( (-1,7,6,-15,-10,10,4,-1)^\top \) & 0.0080532004236864 \\ 
\hline
\( s_{21,0} \) & \( (1,0,0,0,0,0)^\top \) & 0.208277419875994 \\ \hline
\( s_{21,1} \) & \( (-5,0,14,0,-7,1)^\top \) & 0.305356304999981 \\ \hline
\( s_{21,2} \) & \( (9,0,-35,-1,21,-3)^\top \) & 0.0311291821584927 \\ \hline
\( s_{21,3} \) & \( (-5,1,21,1,-14,2)^\top \) & -0.259717694500226 \\ \hline
\( s_{21,4} \) & \( (-6,-5,27,9,-22,3)^\top \) & -0.411902266359988 \\ \hline
\( s_{21,5} \) & \( (13,10,-59,-26,52,-7)^\top \) & -0.344173760232181 \\ \hline
\( s_{21,6} \) & \( (-10,-9,45,26,-45,6)^\top \) & -0.0926921719855542 \\ \hline
\( s_{21,7} \) & \( (1,0,0,0,0,0)^\top \) & 0.208277419875994 \\ \hline
\( s_{21,8} \) & \( (5,9,-31,-26,38,-5)^\top \) & 0.398048476985537 \\ \hline
\( s_{21,9} \) & \( (-4,-10,24,25,-31,4)^\top \) & 0.375302942390668 \\ \hline
\( s_{21,10} \) & \( (1,6,-6,-8,8,-1)^\top \) & 0.152184571859762 \\
\hline
\end{tabular}
\end{table}

\vspace{5mm}
\newpage
\section{Infinite Local Complexity}
\noindent In this section, we apply Danzer's algorithm \cite{Danz02} to prove that every tiling in the tiling space of $\sigma_{n}$, for $n \in \{13, 17, 21\}$, has infinite local complexity.  
The main idea of the algorithm is to select a vertex in a misfit situation and then iterate the substitution along a line passing through this vertex. This procedure yields an infinite sequence of segments of distinct lengths, each corresponding to a unique edge type, and the non-Pisot property of the substitution factor plays a key role in establishing the unbounded growth of distinct lengths.

\begin{thm}
\label{thm:n-foldILC}
Every tiling in $\mathbb{X}_{\sigma_{13}}$ has infinite local complexity.
\end{thm}
\begin{proof}
Consider the edge type $P$ given in Figure~\ref{fig:iteration1}. Such an edge type is legal since a rotated copy of it occurs at the right corner of $\sigma_{13}(T_{13,1})$. Applying the substitution to $P$ produces a vertex in a misfit situation. We embed a copy of $\sigma_{13}(P)$ in the plane as follows (see Figure~\ref{fig:iteration1}): place the chosen misfit vertex at the origin $O$; align the edge of the tile $T$ containing this vertex (but not as one of its corners) along the $x$-axis; and position $T$ entirely within the upper half-plane.

\begin{figure}[H]
\centering
\includegraphics[width=0.8\textwidth]{ILCFig12.png}
			\caption{The image of the edge type $P$ under $\sigma_{13}$, with the misfit vertex placed at the origin and tile $T$ lying above the $x$-axis. (The color opacity has been adjusted to emphasize the segments.)}
\label{fig:iteration1}
\end{figure}

Let us denote by $A_1$ and $B_1$ the endpoints of the edge of the tile $T$ that contains the misfit vertex, as shown in the figure. Applying $\sigma_{13}$ to $\sigma_{13}(P)$ replaces the edge $\overline{A_1B_1}$ with a finite sequence of edges from substituted tiles in the upper half-plane. The total length of these edges is $\mu_{13} |\overline{A_1B_1}|$, where $|\overline{A_1B_1}|$ denotes the length of the edge. Since the origin lies between $A_1$ and $B_1$, one of these new edges must contain the origin. In Figure~\ref{fig:iteration2}, this edge is the base of a tile $T'$, with endpoints $A_2$ and $B_2$.

\begin{figure}[H]
\centering
\includegraphics[width=0.8\textwidth]{ILCFig13.png}
			\caption{Edges $\overline{A_1 B_1}$ and  $\overline{A_2 B_2}$.}
\label{fig:iteration2}
\end{figure}

Iterating the process yields a sequence of edges $(\overline{A_rB_r}){r \in \mathbb{N}}$, where each $\overline{A_rB_r}$ is the base of a tile lying in the upper half-plane and contained in $\sigma_{13}^r(P)$. As illustrated in Figures~\ref{fig:iteration1} and \ref{fig:iteration2}, which show the cases $r = 1$ and $r = 2$, the point $B_r$ is defined so that it lies strictly to the right of the origin. So $|\overline{OB_r}| > 0$ for every $r$.

We now associate to each edge $\overline{OB_r}$ a unique edge type $E_r$, defined as the pair of tiles that intersect along $\overline{OB_r}$ and both contain the origin. For example, $E_1 = \{T, T''\}$ and $E_2 = \{T', T'''\}$ (see Figure~\ref{fig:iteration2}). Since $\sigma_{13}$ is primitive, every edge type $E_r$ appears in every tiling $\mathcal{T} \in \mathbb{X}_{\sigma_{13}}$. Furthermore, if $|\overline{OB_{r_1}}| \neq |\overline{OB_{r_2}}|$ for some $r_1, r_2 \in \mathbb{N}$, then the corresponding edge types $E_{r_1}$ and $E_{r_2}$ are not equivalent. Therefore, if the values $|\overline{OB_r}|$ take infinitely many distinct values, the tiling $\mathcal{T}$ must contain infinitely many non-equivalent edge types, which implies $\mathcal{T}$ exhibits ILC. \\

So, in the following, we aim to show that the values $|\overline{OB_r}|$ take infinitely many distinct values. As a first step, we show that each length $|\overline{OB_r}|$ lies in $\mathbb{Z}[\mu_{13}]$ for all $r \in \mathbb{N}$. \\

We begin with the case $r=1$. We note that $|\overline{OB_1}|=|\overline{A_1B_1}|-|\overline{A_1O}|=s_{13,1}-(s_{13,4}+s_{13,6} )$ (see Figure~\ref{fig:iteration1}; $\overline{A_1B_1}$ is the base of the tile $T$ equivalent to $T_{13,1}$ in the upper half-plane, and $\overline{A_1O}$ is a composition of bases of tiles equivalent to $T_{13,4}$ and $T_{13,6}$ in the lower half-plane). By Lemma~\ref{lem:s_n,j}, we know that $s_{13,j} \in \mathbb{Z}[\mu_{13}]$ for each $j \in \{0,1,\ldots,6\}$. Thus, $|\overline{OB_1}| =s_{13,1}-(s_{13,4}+s_{13,6} ) \in \mathbb{Z}[\mu_{13}]$ by the closure property of $\mathbb{Z}[\mu_{13}]$.  

Now, for $r>1$, the length $|\overline{OB_r}|$ can be expressed using the recursive formula:

\begin{equation}
|\overline{OB_r}|=\mu_{13} |\overline{OB_{r-1}}|-t_r,           
\label{eq:recursion}
\end{equation}
\noindent where $t_r$ is the sum of finitely many prototile edge lengths. (The recursive formula is illustrated in Figure~\ref{fig:proofn=13} for $1< r \leq5$.  So $t_r \in \mathbb{Z}[\mu_{13}]$.) Therefore, since $|\overline{OB_1}|, \ t_r \in \mathbb{Z}[\mu_{13}]$, $|\overline{OB_r}|=\mu_{13}|\overline{OB_{r-1}}|-t_{r} \in \mathbb{Z}[\mu_{13}]$ by the closure property of $\mathbb{Z}[\mu_{13}]$.

\begin{figure}[H]
\centering
\includegraphics[width=0.8\textwidth]{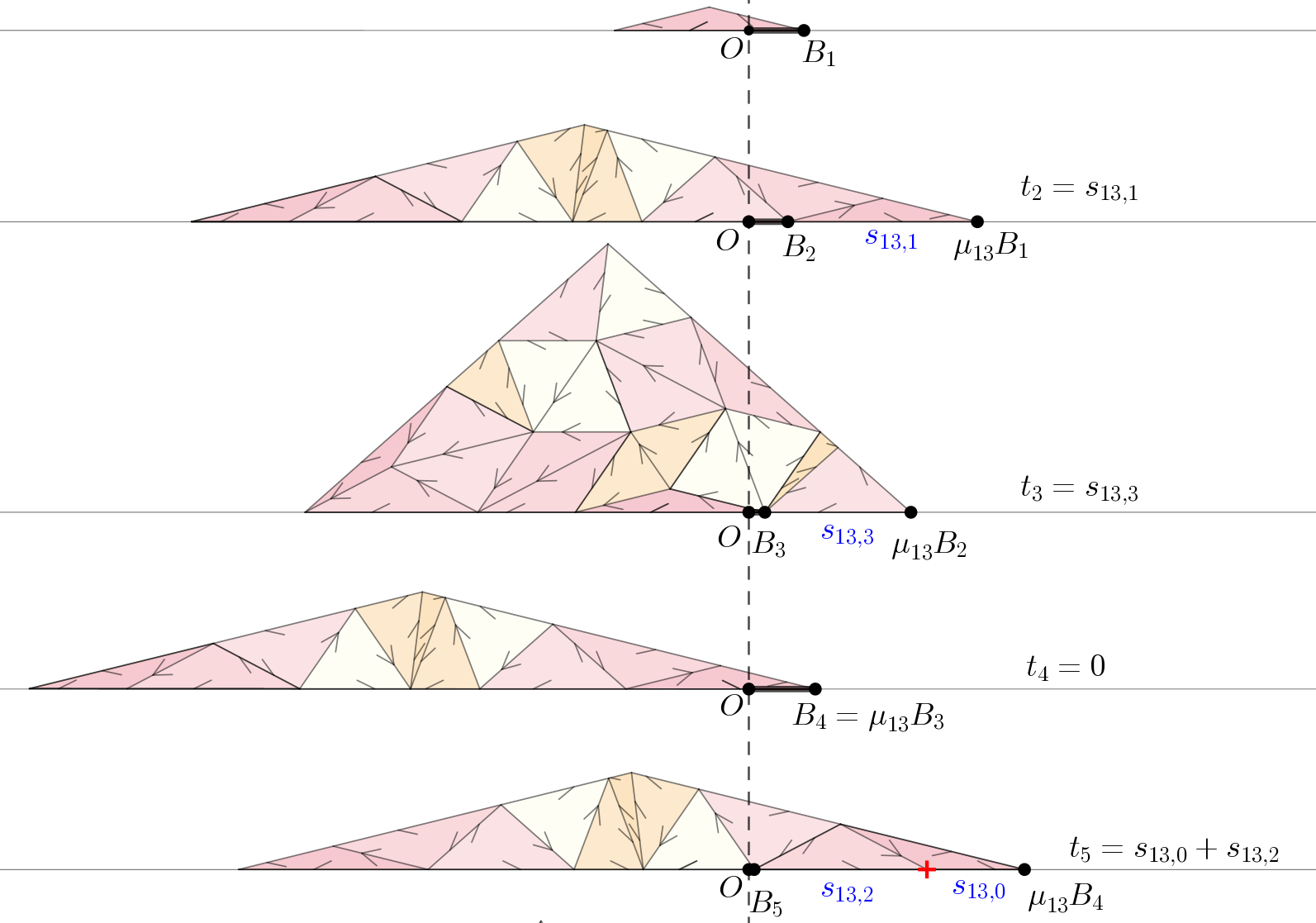}
			\caption{$|\overline{OB_r}|$, $r=1,2,3,4,5.$}
\label{fig:proofn=13}
\end{figure}

Let $w_r = (b_{r_0}, b_{r_1}, \dots, b_{r_5})^{\top} \in \mathbb{Z}^6$ be the coordinate vector of $|\overline{OB_r}|$. As mentioned ealier, multiplication by $\mu_{13}$ in $\mathbb{Z}[\mu_{13}]$ corresponds to left multiplication by the companion matrix $G_{13}$ of the minimal polynomial $g_{13}(x)$ in $\mathbb{Z}^6$. So in $\mathbb{Z}^6$, the recursive formula \ref{eq:recursion} is equivalent to

\begin{equation}
w_r=G_{13} w_{r-1}-u_r,
\label{eq:recursionvec}     	
\end{equation}

\noindent where $u_{r} \in \mathbb{Z}^6$ is the coordinate vector of $t_r$ in $\mathbb{Z}^6$.  To illustrate this, recall that $|\overline{OB_1}|=s_{13,1}-(s_{13,4}+s_{13,6})$ and Figure~\ref{fig:proofn=13} shows that $t_2=s_{13,1}$. From Table~\ref{tab:vectorrep}, we have 

$$
w_1 =y_{13,1}-(y_{13,4}+y_{13,6})=
\begin{bmatrix}
-4 \\ 1 \\ 10 \\ -3 \\ -4 \\ 1
\end{bmatrix}
-
\left(
\begin{bmatrix}
-5 \\ 0 \\ 16 \\ -5 \\ -8 \\ 2
\end{bmatrix}
+
\begin{bmatrix}
-1 \\ 5 \\ 4 \\ -6 \\ -3 \\ 1
\end{bmatrix}
\right)=\begin{bmatrix} 2 \\ -4 \\ -10 \\ 8 \\ 7 \\ -2\end{bmatrix}
 \quad \mbox{and} \quad u_{2} = y_{13,1}= \begin{bmatrix}
-4 \\ 1 \\ 10 \\ -3 \\ -4 \\ 1
\end{bmatrix}.
$$
So, 
\[
w_2 = G_{13} w_1 - u_2 = 
\begin{bmatrix}
0 & 0 & 0 & 0 & 0 & 1 \\
1 & 0 & 0 & 0 & 0 & 1 \\
0 & 1 & 0 & 0 & 0 & -5 \\
0 & 0 & 1 & 0 & 0 & -4 \\
0 & 0 & 0 & 1 & 0 & 6 \\
0 & 0 & 0 & 0 & 1 & 3 
\end{bmatrix}
\begin{bmatrix} 2 \\ -4 \\ -10 \\ 8 \\ 7 \\ -2\end{bmatrix} -
\begin{bmatrix} -4 \\ 1 \\ 10 \\ -3 \\ -4 \\ 1 \end{bmatrix} =
\begin {bmatrix} 2 \\-1 \\ -4 \\ 1 \\ 0 \\ 0 \end{bmatrix}.
\]

Since there is a one-to-one correspondence between the lengths $|\overline{OB_r}|$ and their coordinate vectors $w_r$, we can establish ILC by showing that $w_r$ takes infinitely many distinct values. For this purpose, it is more convenient to work with the diagonal matrix $\widetilde{G}_{13} = A^{-1} G_{13} A$ (see \ref{eq:Gtilde}) rather than $G_{13}$. As described by Danzer in \cite{Danz02}, a key advantage of this approach is that the recursive process can be applied directly to the entries of the vectors.

By left-multiplying the recursion $w_r = G_{13} w_{r-1} - u_r$ by $A^{-1}$ and inserting the identity $I = A A^{-1}$ between $G_{13}$ and $w_{r-1}$, the recursion \eqref{eq:recursionvec} becomes
\begin{equation}
\label{eq:recursionsecond}
w_r' \;=\; \widetilde{G}_{13}\, w_{r-1}' \;-\; u_r',
\end{equation}
where $w_r'= A^{-1} w_r$, $u_r'= A^{-1} u_r$, and $\widetilde{G}_{13}= A^{-1} G_{13} A$.

This is the step where the non-Pisot property of the substitution factor plays a significant role. We focus on the second entry $w_{r,2}'$ of $w_r'$, which is described by
\begin{equation}
\label{eq:recursionsecondeig}
w_{r,2}' = \mu_{13,2} , w_{{r-1},2}' - u_{r,2}',
\end{equation}

\noindent where $u_{r,2}'$ is the second entry of $u_r'$. Recall that $\mu_{13,2} = -\dfrac{1}{s_{13,5}} \approx -1.41002$. We will show that the sequence $(w_{r,2}')_{r\in\mathbb{N}}$ contains infinitely many distinct terms.  

Since there are only a finite number of prototiles, $t_r$ can take only finitely many distinct values. Consequently, the set $\{u_r \mid r \in \mathbb{N}\}$ of coordinate vectors of $t_r$ is also finite, and therefore $\{u_{r,2}' = v_{13,2} u_r \mid r \in \mathbb{N}\}$ must be finite as well. Hence, we can choose a value $U_{13}$ such that
$$
U_{13} \geq \max \{ |u_{r,2}'| \mid r \in \mathbb{N} \}.
$$
To determine $U_{13}$, we use the values $v_{13,2} y_{13,j}$ derived in the previous section and listed in Table~\ref{tab:vectorrep}. Each $u_{r,2}'$ has the form $\sum_{j=0}^6 c_{13,j} \, v_{13,2} y_{13,j}$, where $c_{13,j} \in \{0,1\}$.

Thus, we have
$$
\max \left\{ \left| \sum_{j=0}^{6} c_{13,j} \, v_{13,2} \, y_{13,j} \right| \;\middle|\; j \in I_{13}\right\} 
\geq 
\max \left\{ \left| u_{r,2}' \right| \;\middle|\; r \in \mathbb{N} \right\}.
$$
Using a Python program \cite{PythonManual}, the maximum value of  $\left| \sum_{j=0}^{6} c_{13,j} \, v_{13,2} \, y_{13,j} \right|$
was computed to be approximately $0.997$ (see the Appendix). We therefore fix $U_{13} = 1.2$, which exceeds the maximum.

By applying the triangle inequality to equation~\ref{eq:recursionsecondeig} and using the fact that $U_{13}> |u_{r,2}'|$ for all $r$, we obtain
\begin{align*}
|w_{r,2}'| &\geq |\mu_{13,2}| \cdot |w_{r-1,2}'| - |u_{r,2}'| \\
&= |w_{r-1,2}'| + (|\mu_{13,2}| - 1) |w_{r-1,2}'| - |u_{r,2}'|\\
&> |w_{r-1,2}'| + (|\mu_{13,2}| - 1) |w_{r-1,2}'| - U_{13}. \label{eq:triangle_ineq}
\end{align*}

Thus, $|w_{r,2}'| >|w_{r-1,2}'|+(|\mu_{13,2}| - 1) |w_{r-1,2}'| - U_{13}.$ In particular, 
$|w_{r,2}'|>|w_{r-1,2}'|$ whenever $(|\mu_{13,2}| - 1) |w_{r-1,2}'| - U_{13}>0$, or equivalently, $(|\mu_{13,2}| - 1) |w_{r-1,2}'| >U_{13}$. Since $|\mu_{13,2}|>0$ the condition becomes
\begin{equation}
\label{eq:danzer}
\frac{|w_{r-1,2}'|}{U_{13}} = \frac{|w_{r-1,2}'|}{1.2} > \frac{1}{|\mu_{13,2}| - 1} \approx 2.4389.
\end{equation}

We refer to the inequality \ref{eq:danzer} as \textbf{Danzer’s criterion}. If a positive integer $m$ satisfies the Danzer's inquality, that is,  $\frac{|w_{m-1,2}'|}{1.2} > \frac{1}{|\mu_{13,2}| - 1}\approx 2.4389$, then, as explained, $|w_{m,2}'|>|w_{m-1,2}'|$. It follows that $\dfrac{|w_{m,2}'|}{U_{13}}>\dfrac{|w_{m-1,2}'|}{U_{13}}> \dfrac{1}{|\mu_{13,2}| - 1} \approx 2.4389$, which in turn implies that $|w_{m+1,2}'|>|w_{m,2}'|$. Iterating this argument yields a strictly increasing sequence
$|w_{m-1,2}'| < |w_{m,2}'| < |w_{m+1,2}'| < |w_{m+2,2}'| < \cdots.$ Thus, if we can find such an $m$, the values $|w_{r,2}'|$ take infinitely many values.

Table~\ref{tab:proofn=13} list the values of $|w_{r,2}'|$ and $|w_{r,2}'|/U_{13,2}, r=1,2, \ldots ,8.$ Observe that $|w_{r,5}'|/U_{13,2} \approx 3.2312 >1/ (|\mu_{13,2}|-1) \approx 2.4389$, so $m=5$ satisfies the Danzer's criterion. Therefore, the sequence  $\big(|w_{r,2}'|\big)_{r \in \mathbb{N}}$ has infinitely many distinct terms. Consequently, the sequence $(w_r)_{r \in \mathbb{N}}$ also has infinitely many distinct terms, and so does $\big(|\overline{OB_r}|\big)_{r \in \mathbb{N}}$. Therefore, $\mathcal{T} \in \mathbb{X}_{\sigma_{13}}$ has infinite local complexity.

\renewcommand{\arraystretch}{1.4}
\begin{table}[htbp]
\centering
\caption{$t_r$,$u_r$,$w_r$,$|w_{r,2}'|$ and $\dfrac{|w_{r,2}'|}{U_{13}}$, $r=1,2,\ldots,12.$}
\label{tab:proofn=13}
\resizebox{\textwidth}{!}{%
\scriptsize{
\begin{tabular}{|l|l|l|l|l|l|}
\hline
$r$ & $t_r$ & $u_r$ & $w_r$ & $|w_{r,2}'|$  & $\dfrac{|w_{r,2}'|}{U_{13}}$ \\
&&&$w_r=G_{13} w_{r-1}-u_r, r>1$&& \\ \hline
1 & & & $(2,-4,-10,8,7,-2)^\top$  & 0.744259419 & 0.6202162 \\
2 & $s_{13,1}$ & $(-4,1,10,-3,-4,1)^\top$ & $(2,-1,-4,1,0,0)^\top$ & 1.318017974 & 1.09834831 \\
3 & $s_{13,3}$ & $(0,5,5,-11,-6,2)^\top$ & $(0,-3,-6,7,7,-2)^\top$ & 2.062277393 & 1.7185645 \\
4 & 0 & $(0,0,0,0,0,0)^\top$ & $(-2,-2,7,2,-5,1)^\top$ & 2.907852469 & 2.42321039 \\
\rowcolor{gray!20}
5 & $s_{13,0}+s_{13,2}$ & $(6,-4,-20,11,11,-3)^\top$ & $(-5,3,13,-8,-3,1)^\top$ & 3.877455508 & 3.2312129 \\
6 & 0 & $(0,0,0,0,0,0)^\top$ & $(1,-4,-2,9,-2,0)^\top$ & 5.467290003 & 4.556075 \\
7 & 0 & $(0,0,0,0,0,0)^\top$ & $(0,1,-4,-2,9,-2)^\top$ & 7.708988515 & 6.4241571 \\
8 & $s_{13,0}+s_{13,2}$ & $(6,-4,-20,11,11,-3)^\top$ & $(-8,2,31,-7,-25,6)^\top$ & 11.09250313 & 9.24375261 \\ \hline
\end{tabular}}}
\end{table}
\end{proof}

\begin{rem}
\label{rem:condition}
We summarize below the key conditions used in the proof to establish ILC. 
Certain conditions and steps from \cite{Danz02} were not required in our argument; these will be discussed in the next section.

\begin{enumerate}
\item[(A)] \textbf{Primitivity of the substitution rule.}  
This ensures that the edge types corresponding to the segments $\overline{OB_r}$ occur as legal patches in any tiling associated with $\sigma_{13}$.

\item[(B)] \textbf{Non-Pisot property of the substitution factor.}  
This guarantees the existence of an algebraic conjugate $\mu_{13,2}$ of $\mu_{13}$ with modulus greater than $1$, a key hypothesis for applying Danzer’s criterion.

\item[(C)] \textbf{Distinct roots of the minimal polynomial of the substitution factor.}  
This implies that the companion matrix $G_{13}$ is diagonalizable. Consequently, the recursion for $w_r$ in \eqref{eq:recursionvec} can be transformed into the recursion for $w_r'$ in \eqref{eq:recursionsecond}, allowing the iteration to be analyzed via the second entry $w_{r,2}'$ in \eqref{eq:recursionsecondeig}.

\item[(D)] \textbf{Edge lengths in $\mathbb{Z}[\mu_{13}]$}  (see (E)(a)).

\item[(E)] \textbf{Existence of a misfit vertex in a legal patch} satisfying the following:
\begin{enumerate}
\item[(a)] The length $|\overline{OB_1}|$ lies in $\mathbb{Z}[\mu_{13}]$. Combined with (D), this ensures that $|\overline{OB_r}| \in \mathbb{Z}[\mu_{13}]$ for all $r$.
\item[(b)] Danzer’s criterion is satisfied for some $r$. This yields a monotonic sequence of lengths, which in turn implies that there are infinitely many edge types up to equivalence.
\end{enumerate}
\end{enumerate}
\end{rem}

The steps in Theorem~\ref{thm:n-foldILC} apply analogously to the cases \(n = 17\) and \(n = 21\). To avoid repeating the full argument, we summarize the steps and incorporate Remark~\ref{rem:condition}. Throughout, we use analogous notation, but with \(n\) appended to each symbol; for example, we write \(\overline{OB_{n,r}}\) in place of \(\overline{OB_{r}}\) from Theorem~\ref{thm:n-foldILC}.

\begin{thm}
Every tiling in $\mathbb{X}_{\sigma_{n}}$, $n \in \{17,21\}$, has infinite local complexity.
\end{thm}

\noindent\textit{Sketch of the proof.}
Let \( n \in \{17, 21\} \). 

\begin{figure}[H]
\raggedright
\includegraphics[width=0.8\textwidth]{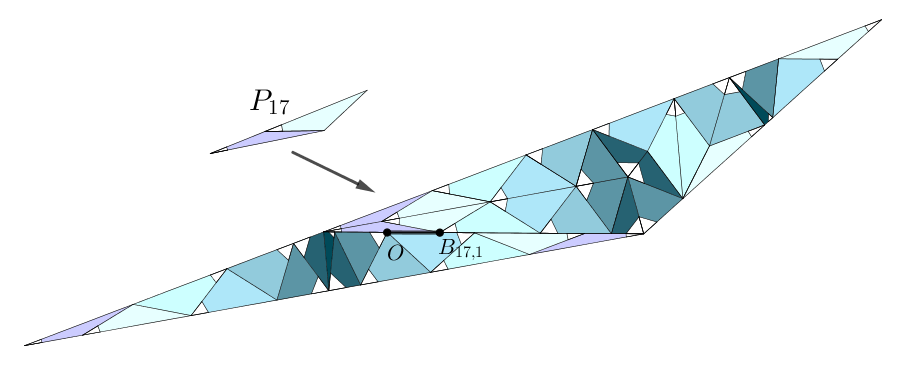}\\
\raggedright
\includegraphics[width=1\textwidth]{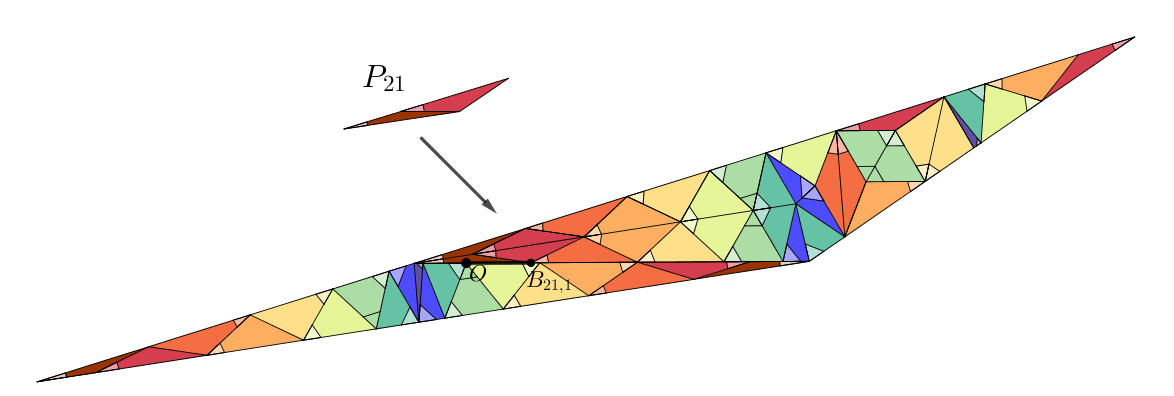}
			\caption{The segments $\overline{OB_{n,1}}, n\in \{17,21\}$.}
\label{fig:misfitn=1721}
\end{figure}

\begin{enumerate}[leftmargin=20mm]
\item[Step 1.] Check conditions (A) to (C).\\[2mm]
The substitution \(\sigma_n\) is primitive since every prototile appears in each 1-order supertile (see Figures~\ref{fig:17-fold} and \ref{fig:21-fold}). Moreover, by Lemma~\ref{lem:mu_n_properties}, the substitution factor \(\mu_n\) is non-Pisot, and the minimal polynomial of \(\mu_n\) has no repeated roots.

\item[Step 2.] Check condition (D) and derive the coordinate vector $y_{n,j}$ of each $s_{n,j}$.\\[2mm]
By Lemma~\ref{lem:s_n,j}, each prototile edge length \(s_{n,j}\) belongs to \(\mathbb{Z}[\mu_n]\). The coordinate vectors $y_{n,j}$ are given in Table~\ref{tab:vectorrep}.\\

\item[Step 3.] Take the algebraic conjugate $\mu_{n,2}$ of $\mu_n$ with modulus greater than $1$, and define
\[
v_{n,2} =
\left( \frac{1}{\mu_{n,2}^{\,d_n - 1}},\,
       \frac{1}{\mu_{n,2}^{\,d_n - 2}},\,
       \ldots,\,
       \frac{1}{\mu_{n,2}},\,
       1 \right).
\]
Compute $v_{n,2}y_{n,j}$ and choose
\[
U_n \;\ge\; \left| \sum_{j=0}^{\frac{n-1}{2}} c_{n,j}\, v_{n,2}\, y_{n,j} \right|
\quad \text{for all } j \in I_n = \{0,1,\ldots,(n-1)/2\} \text{ with } c_{n,j} \in \{0,1\}.
\]\\[2mm]

Using the values of $v_{n,2}y_{n,j}$ given in Table~\ref{tab:vectorrep} and with the aid of a Python program, the computed maximum value of
\[
\left| \sum_{j=0}^{\frac{n-1}{2}} c_{n,j}\, v_{n,2}\, y_{n,j} \right|
\]
is approximately \(0.1066\) for \(n=17\) and \(1.6786\) for \(n=21\). Based on these results, we safely choose \(U_{17} = 0.12\) and \(U_{21} = 1.8\).\\

\item[Step 4.] Choose a vertex in a misfit situation and check whether the corresponding $|\overline{OB_{n,1}}|$ lies in \(\mathbb{Z}[\mu_n]\).\\

Consider Figure~\ref{fig:misfitn=1721}. Applying the substitution to the edge type \(P_n\), which occurs in \(\sigma_n(T_{n,1})\), produces the vertex in a misfit situation indicated in the figure. The segment \(\overline{OB_{n,1}}\) (labeled in the figure) has length
\[
|\overline{OB_{n,1}}|
= s_{n,1} - \bigl(s_{n,\frac{n-5}{2}} + s_{n,\frac{n-1}{2}}\bigr),
\]
which lies in \(\mathbb{Z}[\mu_n]\), since each \(s_{n,j}\) is an element of \(\mathbb{Z}[\mu_n]\).

\item[Step 5.] Iterate the substitution on \(\overline{OB_{n,1}}\) to obtain the values \(\overline{OB_{n,r}}\), $t_{n,r}$, $u_{n,r}$, $w_{n,r} = G_n w_{n,r-1} - u_{n,r}$, and $w_{n,r,2}' = v_{n,2} w_{n,r}$, and then search for \(r\) that satisfies Danzer's criterion:
\begin{equation}
\label{eq:danzer2}
\frac{|w_{n, r, 2}'|}{U_{n}}  > \frac{1}{|\mu_{n,2}| - 1}.
\end{equation}

The values of $t_{n,r}$ are shown in Figures~\ref{fig:n=17nonflc} and~\ref{fig:n=21nonflc}. The companion matrices of $g_{17}(x)$ and $g_{21}(x)$ are given below:
\[
G_{17}=
\begin{bmatrix}
0 & 0 & 0 & 0 & 0 & 0 & 0 & -1\\
1 & 0 & 0 & 0 & 0 & 0 & 0 & -1\\
0 & 1 & 0 & 0 & 0 & 0 & 0 & 7\\
0 & 0 & 1 & 0 & 0 & 0 & 0 & 6\\
0 & 0 & 0 & 1 & 0 & 0 & 0 & -15\\
0 & 0 & 0 & 0 & 1 & 0 & 0 & -10\\
0 & 0 & 0 & 0 & 0 & 1 & 0 & 10\\
0 & 0 & 0 & 0 & 0 & 0 & 1 & 4
\end{bmatrix},
\qquad
G_{21}=
\begin{bmatrix}
0 & 0 & 0 & 0 & 0 & -1\\
1 & 0 & 0 & 0 & 0 & 1\\
0 & 1 & 0 & 0 & 0 & 6\\
0 & 0 & 1 & 0 & 0 & -6\\
0 & 0 & 0 & 1 & 0 & -8\\
0 & 0 & 0 & 0 & 1 & 8
\end{bmatrix}.
\]

We note that \(\frac{1}{|\mu_{17,2}| - 1} \approx 1.20909519021763\) and \(\frac{1}{|\mu_{21,2}| - 1} \approx 2.7131\). As shown in Tables~\ref{tab:proofn=17} and~\ref{tab:proofn=21}, Danzer’s criterion is satisfied for \(\sigma_{17}\) when \(r = 3\), and for \(\sigma_{21}\) when \(r = 9\).

\end{enumerate}
\renewcommand{\arraystretch}{1.4}
\begin{table}[htbp]
\centering
\caption{$u_{17,r}$,$w_{17,r}$,$|w_{17,r,2}'|$ and $\dfrac{|w_{17,r,2}'|}{U_{17}}$, $r=1,2,3$}
\label{tab:proofn=17}
\resizebox{\textwidth}{!}{%
\scriptsize{
\begin{tabular}{|l|l|l|l|l|l|}
\hline
$r$  & $u_{17,r}$ & $w_{17,r}$ & $|w_{17,r,2}'|$  & $\dfrac{|w_{17,r,2}'|}{U_{17}}$ \\ \hline
1 & & $(9, -20, -63, 80, 107, -75, -39, 9)^\top$  & $0.061999291$& $0.5166608$\\
2  & $(-20, 15, 105, -69, -146, 73, 46, -10)^\top$ & $(11, -15, -62, 60, 91, -56, -31, 7)^\top$ & 1.318017974 & 1.1298434 \\
\rowcolor{gray!20}
3  &$(-19, 15, 105, -69, -146, 73, 46, -10)^\top$  &$(12, -11, -71, 49, 101, -52, -32, 7)^\top$ & 0.2106974 & 1.7558120 \\ \hline
\end{tabular}}}
\end{table}

\renewcommand{\arraystretch}{1.4}
\begin{table}[htbp]
\centering
\caption{$u_{21,r}$,$w_{21,r}$,$|w_{21,r,2}'|$ and $\dfrac{|w_{21,r,2}'|}{U_{21}}$, $r=1,2,\ldots,9$}
\label{tab:proofn=21}
\resizebox{\textwidth}{!}{%
\scriptsize{
\begin{tabular}{|l|l|l|l|l|l|}
\hline
$r$ & $u_{21,r}$ & $w_{21,r}$ & $|w_{21,r,2}'|$  & $\dfrac{|w_{21,r,2}'|}{U_{21}}$ \\ \hline
1 & & $(-11, -15, 51, 34, -53, 7)^\top$  & $0.244876744$& $0.136042636$\\
2  & $(0, 1, 0, 0, 0, 0)^\top$  & $(-7, -5, 27, 9, -22, 3)^\top$ & 0.620179686  & 0.34454427\\
3  & $(4, 0, -21, -1, 14, -2)^\top$ & $(-7, -4, 34, 10, -29, 4)^\top$ &1.185253686  &0.65847427 \\
4  & $(13, 10, -59, -26, 52, -7)^\top$ & $(-17, -13, 79, 36, -74, 10)^\top$ &1.277945858  &0.709969921 \\
5  & $(-11, -4, 41, 9, -29, 4)^\top$  & $(1, -3, 6, 10, -15, 2)^\top$  &1.927475923 &1.070819957 \\
6  & $(-6, 2, 21, 1, -14, 2)^\top$ & $(4, 1, -12, -7, 8, -1)^\top$ &2.454963437  &1.363868576 \\
7  & $(-10, 1, 35, 1, -21, 3)^\top$ & $(11, 2, -40, -7, 22, -3)^\top$ &3.405463094  &1.891923941 \\
8  & $(-6, -5, 27, 9, -22, 3)^\top$  & $(9, 13, -43, -31, 39, -5)^\top$  &4.248761151  &2.360422862 \\
\rowcolor{gray!20}
9  & $(-10, 1, 35, 1, -21, 3)^\top$  & $(15, 3, -52, -14, 30, -4)^\top$ &5.860426532  &3.255792518 \\ \hline
\end{tabular}}}
\end{table}
	
\begin{figure}[H]
\centering
\includegraphics[width=0.8\textwidth]{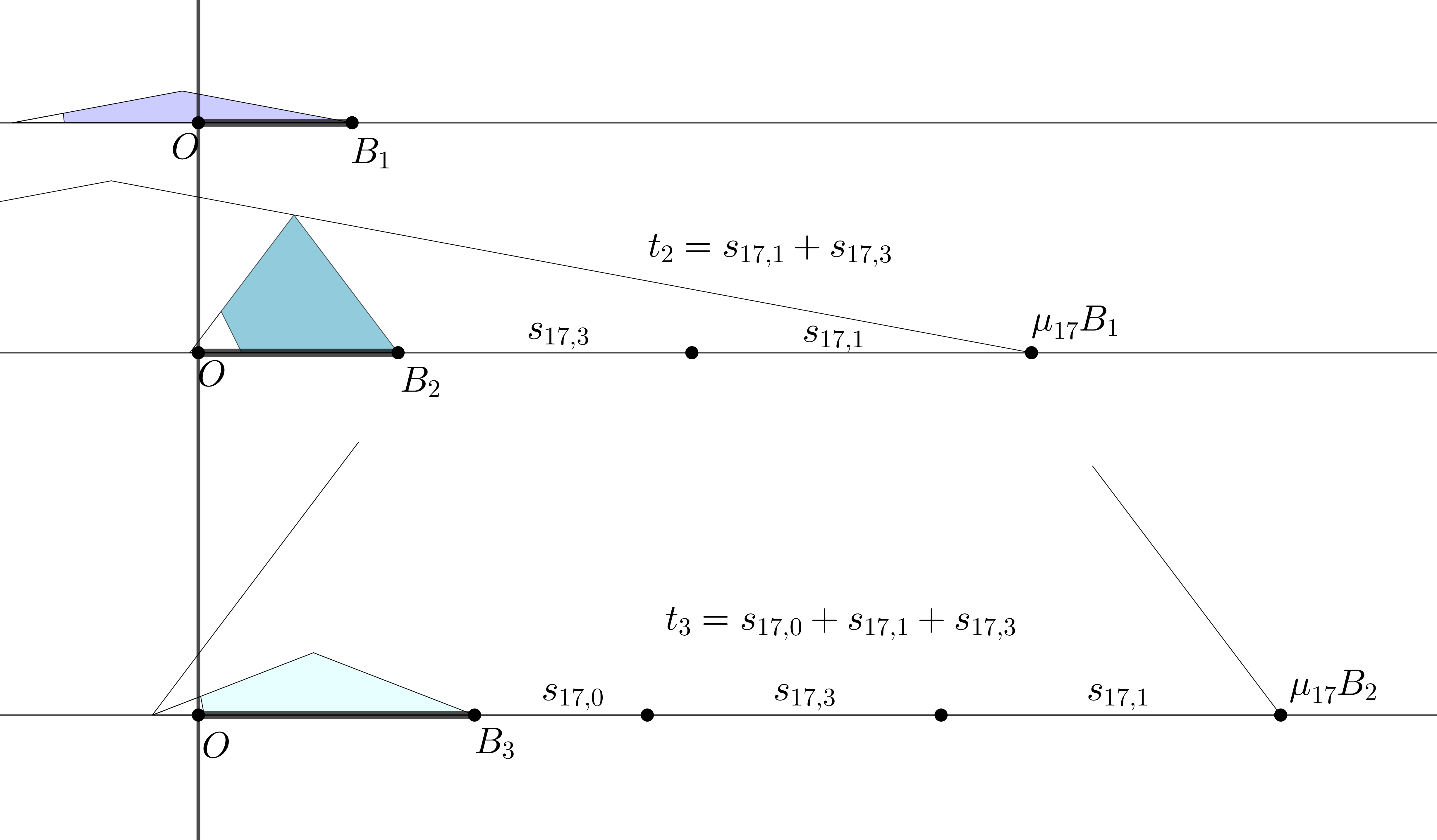}
\caption{The segments $\overline{OB_{17,r,1}}, r=1,2,3$.}
\label{fig:n=17nonflc}
\end{figure}

\begin{figure}[H]
\centering
\includegraphics[width=1\textwidth]{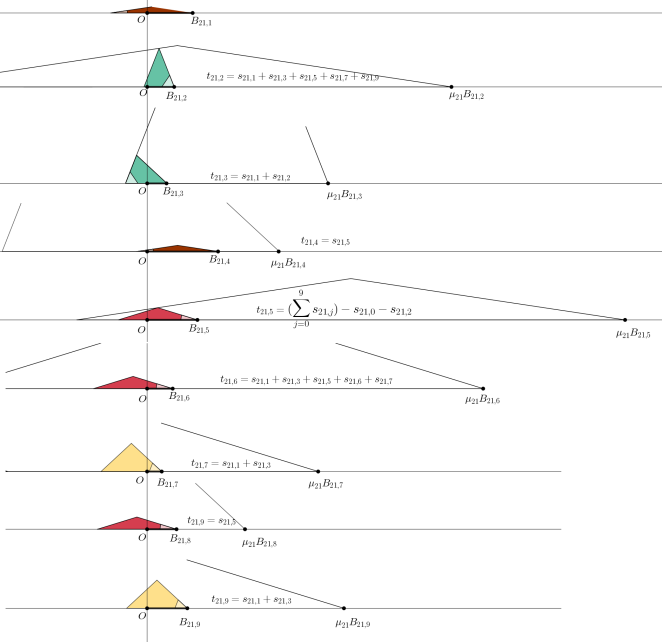}
\caption{The segments $\overline{OB_{21,r,1}}, r=1,2,\ldots,9$.}
\label{fig:n=21nonflc}
\end{figure}

\newpage
\section{Additional Remarks and Conjectures}
\begin{rem}
We describe here some modifications we applied to the algorithm given in \cite{Danz02}.
\begin{enumerate}
\item In \cite{Danz02}, the corresponding constant \( U_{13,2} \) satisfies
\[
U_{13,2} = \max \{ |u_{r,2}'| \mid r \in \mathbb{N} \}.
\]
As shown in the proof, we relaxed this condition by replacing it with the inequality
\[
U_{13,2} \geq \max \{ |u_{r,2}'| \mid r \in \mathbb{N} \}.
\]
This relaxation still achieves the desired result. The main advantage of using the exact value is that it may satisfy Danzer’s criterion~\ref{eq:danzer} for a smaller \( r \). This is evident in Table~\ref{tab:proofn=13}, where the increasing behavior of \( |w_{r,2}'| \) begins from \( r = 1 \). However, determining \( \max \{ |u_{r,2}'| \mid r \in \mathbb{N} \} \) requires iterating the substitution on the edge \( \overline{OB_{1}} \) until all exact values of \( t_r \) are obtained, which is more labor-intensive than simply considering the maximum \textit{possible} value of \( |u_{r,2}'| \).

\item The matrix \(M_1\) in \cite{Danz02}, which was used to obtain the minimal polynomial (and is equivalent to \(M_n\) in this work), is the substitution matrix of the one-dimensional substitution on the set of edges naturally induced by the substitution rule. For convenience, we did not use such a matrix here, since \(M_n\) already contains all the relevant information for the proof, such as the minimal polynomial and the eigenvector containing the lengths of all prototile edges.

\item We did not explicitly use condition (A3) (in our notation, this is $\dfrac{|\overline{A_1B_1}|}{|\overline{OB_1}|} \notin \mathbb{Q}(\mu_n)$) in \cite{Danz02}. This is because each edge length \(s_{n,j}\) (see Lemma~\ref{lem:s_n,j}) and \(|OB_{n,1}|\) lies in \(\mathbb{Z}[\mu_n]\), which in turn implies, by the closure property, that \(|OB_{n,r}| \in \mathbb{Z}[\mu_n]\) for all \(r\). This is precisely the reason why (A3) is included as a condition in \cite{Danz02}, and hence it is not necessary in our proof.
\end{enumerate}
\end{rem}

\begin{figure}[H]
\includegraphics[width=1\textwidth]{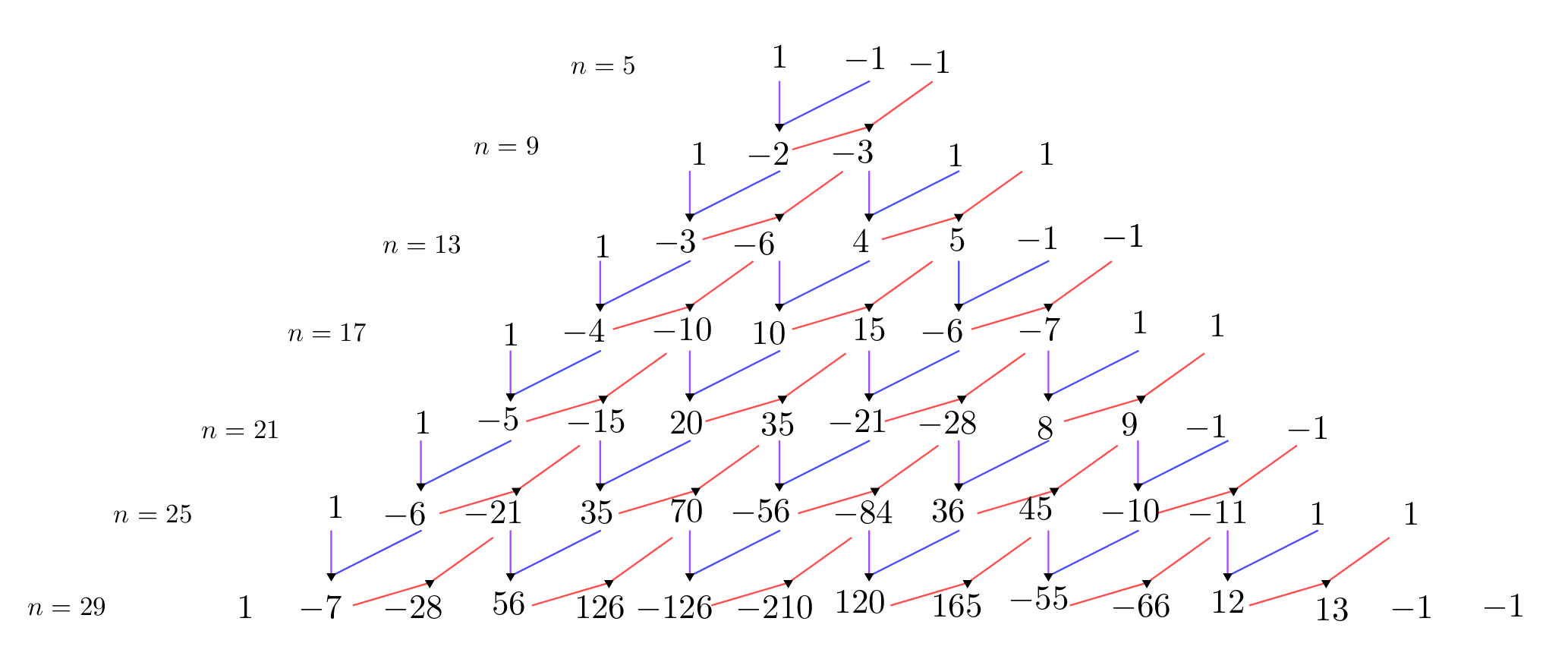}
			\caption{An array of numbers describing the coefficients of  $q_n (x)$ in standard form.}
\label{fig:pattern}
\end{figure}

\begin{conj}
\label{thm:Conjecture1}
During the calculation of the characteristic polynomial $p_n (x)$, a pattern in its coefficients becomes apparent. This polynomial can be expressed as $p_n (x)=-(x-1) q_n (x)$, where the coefficients of $q_n (x)$ in standard form follow the pattern described in Figure~\ref{fig:pattern}, at least for the given values of $n$. The signs of the coefficients from left to right follow a specific sequence: $+,-,-,+,+,-,-,+,+,\ldots$ The last two coefficients in each row are equal to either $1$ or $-1$, depending on $n$. The remaining coefficients are obtained by summing the absolute values of the two numbers indicated by the adjacent segments of the same color, then assigning the appropriate signs according to their positions. Clearly, the minimal polynomial $g_n (x)$ is equal to $q_n (x)$ if $q_n (x)$ is irreducible. An analogous pattern exists for the cases $n\in \{7,11,15,19,23,27,31\}$, and it is the author’s conjecture that this pattern holds for all odd $n\geq 5.$
\end{conj}

\begin{conj}
\label{thm:Conjecture2}
In addition to the previous observation and conjecture, we notice that the roots of the characteristic $p_n (x)$ of $M_n$, at least for the values of $n$ shown in Figure~\ref{fig:pattern}, are given by $\mu_n=\dfrac{1}{s_{n,(n-1)/2}}, -\dfrac{1}{s_{n,(n-3)/2}}, \ldots, \dfrac{1}{s_{n,2}}, -\dfrac{1}{s_{n,1}}, \dfrac{1}{s_{n,0}}=1$. This indicates that each eigenvalue of $M_n$ takes the form $\pm \dfrac{1}{s_{n,j}}$. We conjecture that this holds for all odd $n\geq5$. If this conjecture is correct, then $\mu_n$ is non-Pisot if  $g_n (x)$ is irreducible for odd $n\geq11$, since $\dfrac{1}{s_{n,(n-3)/2}}= \dfrac{1}{2\sin(3\pi/2n)}>1$. 
\end{conj}

\section{Conclusion}
In this work, we have presented explicit examples of substitution tilings with infinite local complexity (ILC) that also exhibit $n$-fold rotational symmetry for $n \in \{13, 17, 21\}$. The proof of ILC relies on Danzer’s algorithm. However, the implementation of this algorithm is subtle and does not automatically yield ILC.

Among the conditions summarized in Remark~\ref{rem:condition}, the most decisive one is the choice of a vertex in a misfit situation. Indeed, some misfit vertices lead to only a finite number of values of the lengths $|\overline{OB_{n,r}}|$. This occurs when, for some $m$, the segment $\overline{OB_{n,m}}$ coincides with one of the prototile edges, or equivalently, when the point $A_m$ coincides with the origin. From that point on, all later values $|\overline{OB_{n,r}}|$ with $r \ge m$ are equal to prototile edge lengths $s_{n,j}$. Since there are only finitely many prototiles, the quantity $|\overline{OB_{n,r}}|$ then takes only finitely many values, and ILC cannot be obtained from that particular misfit vertex. In such a case, one must search for a different misfit situation --- or, as Danzer put it, ``perhaps, another misfit situation will do \textit{better}''.

For the cases $n = 13$ and $n = 17$, all of the misfit vertices considered by the author satisfied Danzer’s criterion after only a few iterations. The situation is different for $n = 21$. As Table~\ref{tab:proofn=21} shows, it takes nine iterations before Danzer’s criterion is satisfied. The main difference lies in the value of $\mu_{n,2}$: while for $n = 13$ and $n = 17$ one has $\mu_{n,2} = -\dfrac{1}{s_{n,(n-3)/2}}$, in the case $n = 21$ this value is a root of the characteristic polynomial of $M_n$, but not of the minimal polynomial of $\mu_n$. We are therefore forced to use $\mu_{21,2} = -\dfrac{1}{s_{21,(21-5)/2}}$, which in turn requires more iterations before Danzer’s criterion is satisfied.

Finally, in light of Conjectures~\ref{thm:Conjecture1} and~\ref{thm:Conjecture2}, and noting that the same uniform misfit vertex is used in the cases considered here (with
$|\overline{OB_{n,1}}| = s_{n,1} - \bigl(s_{n,\frac{n-5}{2}} + s_{n,\frac{n-1}{2}}\bigr)$,
it appears reasonable to expect that any substitution tiling with prototiles $T_{n,1}, T_{n,2}, \ldots, T_{n,(n-1)/2}$ and substitution factor
\[
\mu_n = \frac{1}{2\sin(\pi/2n)}, \qquad n \in \{13,17,21,25,\ldots\},
\]
defined in the same way as the substitutions $\sigma_n$ for $n = 13, 17, 21$, will exhibit infinite local complexity, provided that $\mu_n$ is non-Pisot.  
This claim, however, remains open and requires further study.

\section*{Appendix}

\begin{figure}[H]
\centering
\includegraphics[scale=1]{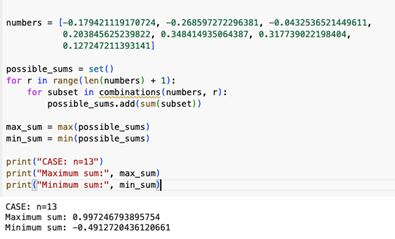}\\
\includegraphics[scale=1]{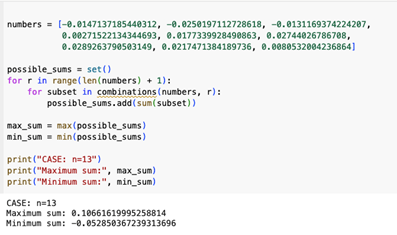}\\
\includegraphics[scale=1]{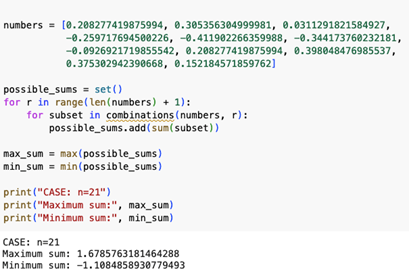}\\
\end{figure}

\end{document}